\documentclass[a4paper]{amsart}
\def\today{April 07, 2006}
\date{\today}
\dedicatory{
 Dedicated to Professor Mitsuhiro Itoh 
 on the occasion of his sixtieth birthday
}
\usepackage[dvips]{graphicx}
\usepackage{amsmath}
\usepackage{verbatim,enumerate}
\usepackage[usenames]{color}
\usepackage{amsthm} 
\theoremstyle{plain}
 \newtheorem{theorem}{Theorem}[section]
 \newtheorem*{theorem*}{Theorem}
 \newtheorem*{proposition*}{Proposition}
 \newtheorem{proposition}[theorem]{Proposition}
 
 \newtheorem{corollary}[theorem]{Corollary}
 \newtheorem{fact}[theorem]{Fact}
 \newtheorem{introtheorem}{Theorem}
 
\theoremstyle{remark}
 \newtheorem{definition}[theorem]{Definition}
 \newtheorem{remark}[theorem]{Remark}
 \newtheorem*{remark*}{Remark}
 \newtheorem{example}[theorem]{Example}
 \newtheorem*{problem}{Problem}
 \newtheorem*{acknowledgement}{Acknowledgement}
\numberwithin{equation}{section}

\newcommand{\Lor}{\boldsymbol{L}}
\newcommand{\Z}{\boldsymbol{Z}}
\newcommand{\Q}{\boldsymbol{Q}}
\newcommand{\R}{\boldsymbol{R}}
\newcommand{\C}{\boldsymbol{C}}
\newcommand{\E}{\mathcal{E}}

\newcommand{\iii}{i\,}
\newcommand{\id}{\operatorname{id}}
\newcommand{\trace}{\operatorname{trace}}
\newcommand{\ord}{\operatorname{ord}}
\newcommand{\SL}{\operatorname{SL}}
\newcommand{\SU}{\operatorname{SU}}
\newcommand{\U}{\operatorname{U}}

\newcommand{\PSL}{\operatorname{PSL}}

\newcommand{\Herm}{\operatorname{Herm}}
\newcommand{\inner}[2]{\left\langle{#1},{#2}\right\rangle}
\newcommand{\Qc}{{Q_\mathrm{c}}}
\newcommand{\Cf}{{C_f}}
\newcommand{\Gc}{{G_\mathrm{c}}}
\newcommand{\G}{\mathcal{G}}
\newcommand{\omegac}{{\omega_\mathrm{c}}}
\newcommand{\thetac}{{\theta_\mathrm{c}}}
\newcommand{\rhoc}{\rho_\mathrm{c}}
\renewcommand{\Mc}{M_{\mathrm{c}}}
\newcommand{\McUniv}{\widetilde M_{\mathrm{c}}}
\newcommand{\Ec}{\E_{\mathrm{c}}}

\newcommand{\orig}{\mathrm{orig}}
\newcommand{\diag}{\operatorname{diag}}
\renewcommand{\Re}{\operatorname{Re}}
\renewcommand{\Im}{\operatorname{Im}}
\newcommand{\Ker}{\operatorname{Ker}}
\title[Flat Fronts]{
   Flat fronts in hyperbolic 3-space\\
   and their caustics
}
\author[M.~Kokubu]{M. Kokubu}
\address[Masatoshi Kokubu]{%
   Department of Natural Science,
   School of Engineering,
   Tokyo Denki University,
   2-2 Kanda-Nishiki-Cho,
   Chiyoda-Ku, Tokyo, 101-8457,
   Japan
}
\email{kokubu@cck.dendai.ac.jp}
\author[W.~Rossman]{W. Rossman}
\address[Wayne Rossman]{%
   Department of Mathematics, Faculty of Science,
   Kobe University,
   Rokko, Kobe 657-8501, Japan%
}
\email{wayne@math.kobe-u.ac.jp}
\author[M.~Umehara]{M. Umehara}
\address[Masaaki Umehara]{%
   Department of Mathematics, Graduate School of Science,
   Osaka University,
   Toyonaka, Osaka 560-0043,
   Japan
}
\email{umehara@math.wani.osaka-u.ac.jp}
\author[K.~Yamada]{K. Yamada}
\address[Kotaro Yamada]{%
   Faculty of Mathematics,
   Kyushu University, 
   Higashi-ku, Fukuoka 812-8581, Japan%
}
\email{kotaro@math.kyushu-u.ac.jp}
\begin{document}
\begin{abstract}
After G\'alvez, Mart\'\i{}nez and Mil\'an 
discovered a (Weierstrass-type) holomorphic representation
formula for flat surfaces in hyperbolic $3$-space $H^3$, 
the first, third and fourth authors here 
gave a framework 
for complete flat fronts with singularities in $H^3$.  In the present 
work we broaden the notion of completeness to {\it weak 
completeness}, and of front to {\it p-front}.
As a front is a p-front and completeness implies weak completeness, 
the new framework and results here apply to a more general class 
of flat surfaces.  

This more general class contains the caustics of flat 
fronts --- shown also to be flat by Roitman (who gave a 
holomorphic representation formula for them) ---
which are an important class of surfaces and are generally 
not complete but only weakly complete.  Furthermore, 
although flat fronts have globally defined normals, caustics 
might not, making them flat fronts only locally, and hence only p-fronts.  
Using the new framework, we obtain characterizations for caustics. 
\end{abstract}
\maketitle
\section{Introduction}
For an arbitrary Riemannian $3$-manifold $N^3$, a $C^\infty$-map 
\[
  f\colon{}M^2\longrightarrow N^3
\]
from a $2$-manifold $M^2$ is a ({\em wave}) {\em front\/}
if $f$ lifts to a smooth immersed section 
\[
  L_f\colon{}M^2\longrightarrow T_1N^3 (\approx T_1^*N^3)
\]
of the unit tangent vector bundle $T_1N^3$
such that $df(X)$ is perpendicular to $L_f(p)$ for all
$X\in T_pM^2$ and $p\in M^2$.  
Fronts generalize immersions, as they allow for singularities. 
The lift $L_f$ can be viewed as a globally defined unit normal 
vector field of $f$.
However, global definedness of $L_f$ can be a stronger 
condition than desired. 
Sometimes p-fronts are more appropriate: The map
$f$ is called a {\em p-front\/} if for each $p\in M^2$,
there is a neighborhood $U$ of $p$ such that the restriction $f|_U$ 
is a front. 
The projectified cotangent bundle $P(T^*N^3)$ has a canonical
contact structure, and a p-front can be considered as the projection 
of a Legendrian immersion of $M^2$ into $P(T^*N^3)$.
A p-front $f$ is a front if and only if there exists a globally
defined unit normal vector field, in which case 
we say $f$ is  {\it co-orientable\/}.  Otherwise, $f$ is
{\it non-co-orientable}.

\begin{figure}
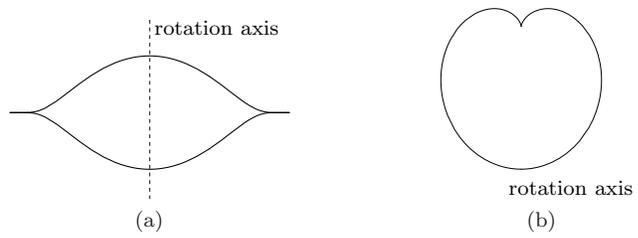

\footnotesize
\begin{center}
   \begin{tabular}{c@{\hspace{2cm}}c}
        \input{figure-1a.pstex_t} &
        \input{figure-1b.pstex_t} \\
        (a) & (b)
  \end{tabular}
\end{center}
\caption{A flying saucer front and toroidal p-front.}\label{fig:curves}
\end{figure}

For example, Figure~\ref{fig:curves} (a) is a plane curve front
with a globally defined  unit normal vector field. 
On the other hand, Figure~\ref{fig:curves} (b) is the lima\c{c}on,
whose unit normal vector field is not single-valued on the curve.
Rotating the first curve about its central vertical
axis gives a surface like a ``flying saucer'', 
which is a front in Euclidean space $\R^3$.
If we rotate the lima\c{c}on 
about an axis disjoint
from it, we get a torus with one cuspidal edge, 
which is a non-co-orientable p-front.

Now let $N^3$ be the hyperbolic 
3-space $H^3$ of constant curvature $-1$, and 
\[
  f\colon{}M^2\longrightarrow H^3(\subset \Lor^4)
\]
a front with a unit normal vector field $\nu\colon{}M^2\to S^3_1$,
where $S^3_1$ denotes the de Sitter space in the Minkowski $4$-space
$\Lor^4$.
Then for each real number $t$ 
\begin{equation}\label{eq:parallel}
\begin{split}
  f_t& =(\cosh t) f + (\sinh t) \nu \\
  \nu_t& =(\cosh t) \nu + (\sinh t) f
\end{split}
\end{equation}
gives a new front called a {\em parallel front\/} of $f$.
If $f$ is a flat surface, then the parallel surfaces $f_t$ are flat as 
well (basic properties of flat surfaces in $H^3$ are in \cite{GMM}, 
\cite{KUY1}, \cite{KUY2}).  
So we say that $f$ is a {\em flat front\/} if for each $p\in M^2$,
there exist a real number $t \in \R$ and a neighborhood $U$ of $p$ so that
the restriction $f_t|_U$ is an immersion with vanishing
Gaussian curvature.  
Hence each parallel front $f_t$ of $f$ 
is also a flat front. 
Moreover, for each non-umbilic point $p\in M^2$, 
there is a unique $t(p) \in \R$ so that
$f_{t(p)}$ is not an immersion at $p$
(see Remark~\ref{rmk:add} for details). 
Then the singular locus (or equivalently, the set of focal points) 
is the image of the map
\[
  C_f\colon{}M^2\setminus\{\text{umbilic points}\}\ni
         p\longmapsto f_{t(p)}\in H^3,
\]
which is called the {\it caustic\/} (or {\it focal surface}) 
of $f$.
Roitman \cite{R} pointed out that $C_f$ is a flat p-front,
and gave a holomorphic representation formula for such caustics. 
We remark that 
{\em caustics  of flat surfaces in $\R^3$ and the $3$-sphere
$S^3$ are also flat.}
However, the caustics of complete flat fronts are not fronts in general,
as the unit normal vector field of $C_f$ might not extend globally.
Moreover, they might not be complete, since the singular set
may accumulate at the ends; they instead satisfy the weaker condition 
{\it weak completeness}. 
The purpose of this paper is to give
a broader framework for flat surfaces in $H^3$ that contains the caustics and 
gives characterizations of them.

After giving preliminaries in Section~\ref{sec:prelim}, 
in Section~\ref{sec:complete} 
we define the notion ``weak completeness'' of fronts.
There we also define $f$ to be of {\it finite type\/} if 
the hermitian part of the first fundamental form
with respect 
to the complex structure induced by the second fundamental form
has finite total curvature, and then prove: 
\begin{introtheorem}\label{thm:A}
 A complete flat front is weakly complete and of finite type. 
 Conversely, if $f\colon{}M^2\to H^3$ is a weakly complete flat 
 front of finite type, then there exists a finite set of real numbers
 $t_1,\dots,t_n$ such that $f_t\colon{}M^2\to H^3$
 is a complete flat front for all $t \in \R \setminus \{t_1,\dots,t_n\}$.
\end{introtheorem}
Section~\ref{sec:p-front} is a study of p-fronts, where we 
prove that any non-co-orientable p-front is the projection of a 
doubly-covering front, and we prove Theorem \ref{thm:B}.  For a regular surface, 
orientability and co-orientability are the same notion, but this is not so for 
p-fronts, as this theorem shows.  
\begin{introtheorem}\label{thm:B}
 Any flat $p$-front is orientable.
\end{introtheorem}
This is an important property of flat surfaces in $H^3$, because, 
there do in fact exist flat M\"obius bands in $\R^3$ and $S^3$. 
(For $S^3$ this is a deep fact, 
since such a front in $S^3$ can be of class $C^\infty$, 
but is never $C^\omega$, see G\'alvez and Mira \cite{GM}.) 

Section~\ref{sec:caustic} summarizes properties of caustics.
In Section~\ref{sec:end}, we investigate ends of 
the caustic $C_f$ of a given flat front $f$.
The ends of $C_f$ come from the umbilic 
points or the ends of $f$, called {\em $U$-ends\/} and {\em $E$-ends} of $C_f$, 
respectively.  Calling an end 
{\it regular\/} if the two hyperbolic Gauss maps have at most poles at 
the end, we prove: 

\begin{introtheorem}\label{thm:C}
 For a non-totally-umbilic
 flat front $f \colon M^2 \to H^3$, the following 
 assertions are equivalent{\rm :}
 \begin{enumerate}
 \item $M^2$ is biholomorphic to 
       $\overline{M}^2 \setminus \{p_1,\dots,p_n \}$ for some compact Riemann 
       surface $\overline{M}^2$ containing the points $p_j$,  
       and $f$ is a weakly complete flat front, 
       all of whose ends are regular. 
 \item The caustic $C_f$ is a weakly complete p-front
       of finite type, all of whose ends are regular.  
 \end{enumerate}
\end{introtheorem}

\begin{remark*}
1) The asymptotic behavior of weakly complete regular ends will be
   treated in the forthcoming \cite{KRUY}. 
2) Generic singularities of flat fronts in $H^3$ consist of cuspidal
   edges and swallowtails \cite{KRSUY}.
   But although cone-like singularities of fronts are 
   not generic, they are still important.  
   Several remarkable results on flat surfaces with cone-like singularities
   were recently given by G\'alvez and Mira \cite{GM2}.
3) A differential geometric viewpoint of fronts was given in 
   \cite{SUY}, 
   where ``singular curvature'' on cuspidal edges was introduced.
   Cuspidal edges on flat fronts in $H^3$ have negative singular 
   curvature \cite[Theorem~3.1]{SUY}.
\end{remark*}

We also provide new examples, in addition to known examples,
showing that the results here are not vacuous.
We characterize the known flat fronts of revolution and peach fronts in
terms of their caustics, in Section~\ref{sec:caustic}.
In Section~\ref{sec:example}, we prove the general existence of complete 
flat fronts with given orders of ends on arbitrary Riemann surfaces of finite
topology, and in particular give explicit data for examples 
of genus $k$ with $4k+1$ embedded ends for all $k \geq 1$.  
Also, in Section~\ref{sec:p-front}, we give an example of a weakly
complete p-front that is not the caustic of any flat front.

\begin{acknowledgement}
 The third and fourth authors thank A. J. G\'alvez, 
 A. Mart\'\i{}nez and P. Mira for fruitful discussions 
 at Granada University.
 The authors also thank the referee for a very careful reading.
\end{acknowledgement}

\section{Preliminaries}
\label{sec:prelim}
Let $\Lor^4$ be the Minkowski $4$-space with
the inner product $\inner{~}{~}$ of signature $(-,+,+,+)$.
The hyperbolic $3$-space $H^3$ is considered as the upper-half 
component of the two-sheeted hyperboloid in $\Lor^4$
with metric induced by $\inner{~}{~}$.
Identifying $\Lor^4$ with the set of $2\times 2$-hermitian
matrices $\Herm(2)$, we have 
\[
      \Lor^4\ni (x_0,x_1,x_2,x_3)\leftrightarrow
      \begin{pmatrix}
        x_0+x_3 & x_1+i x_2 \\
        x_1-i x_2 & x_0-x_3
      \end{pmatrix}\in\Herm(2),
\]
\begin{align*}
     H^3 &= \{x=(x_0,x_1,x_2,x_3)\in\Lor^4\,;\,
             \inner{x}{x}=-1,x_0>0\}\\
         &= \{ X\in\Herm(2)\,;\,\det X=1,~\trace X>0\}\\
         &= \{ aa^*\,;\,a\in\SL(2,\C)\}=\SL(2,\C)/\SU(2),
\end{align*}
where $a^*={}^t\bar a$.
The Lie group $\SL(2,\C)$ acts isometrically on $H^3$ via 
\begin{equation}\label{eq:isometric}
 X \longmapsto aXa^*\qquad
  \bigl(a\in\SL(2,\C), X\in H^3\subset\Herm(2)\bigr).
\end{equation}
In fact, the identity component of the isometry group of $H^3$
can be identified with $\PSL(2,\C)=\SL(2,\C)/\{\pm\id\}$.

Let $M^2$ be an oriented $2$-manifold, and let 
\[
    f \colon M^2\longrightarrow H^3=\SL(2,\C)/\SU(2)
\]
be a front with Legendrian lift (see \cite{KUY2})
\[
    L_f \colon M^2\longrightarrow T_1^* H^3=\SL(2,\C)/\U(1).
\]
Identifying $T_1^*H^3$ with $T_1H^3$, we can write $L_f=(f,\nu)$,
where $\nu(p)$ is a unit vector in $T_{f(p)}H^3$ such that
$\inner{df(p)}{\nu(p)}=0$ for each $p\in M^2$.
We call $\nu$ a {\em unit normal vector field\/} of the front $f$.

Suppose that $f$ is flat, then there is a (unique) complex structure on
$M^2$, called the {\it canonical complex structure}, that is 
(conformally) compatible 
with the second fundamental form wherever it is definite,
and a holomorphic Legendrian immersion
\begin{equation}\label{eq:lift-eq}
  \E_f\colon \widetilde{M}^2 \longrightarrow \SL(2,\C)
\end{equation}
such that $f$ and $L_f$ are projections of $\E_f$, where 
$\pi\colon\widetilde{M}^2\to M^2$ is the universal cover of
$M^2$.
Here, holomorphic Legendrian map means that $\E_f^{-1}d\E_f$ is 
off-diagonal (see \cite{GMM}, \cite{KUY1}, \cite{KUY2}).  
The map $f$ and its unit normal vector field $\nu$ are 
\begin{equation}\label{eq:front-normal}
   f = \E^{}_f\E^*_f,\qquad
   \nu = \E^{}_fe_3
         \E^*_f,\qquad
   e_3=         \begin{pmatrix}
                 1 & \hphantom{-}0 \\
                 0 & -1
         \end{pmatrix}         .
\end{equation}
If we set
\begin{equation}\label{eq:ode}
 \E_f^{-1}d\E_f=
 \begin{pmatrix}
  0 & \theta \\
  \omega & 0
 \end{pmatrix},
\end{equation}
the first and second fundamental forms $ds^2$ and $dh^2$ are given by 
\begin{equation}\label{fff-sff}
\begin{array}{rl}
 ds^2&=|\omega+\bar \theta|^2= Q + \overline{Q} + 
  (|\omega|^2+|\theta|^2), \qquad Q=\omega\theta \; , \\
 dh^2&=|\theta|^2-|\omega|^2
\end{array}
\end{equation}
for holomorphic $1$-forms $\omega$ and $\theta$ defined on 
$\widetilde{M}^2$, with 
$|\omega|$ and $|\theta|$ defined on $M^2$ itself.  
We call $\omega$ and $\theta$ the {\em canonical forms\/} of the 
front $f$ (or the Legendrian immersion $\E^{}_f$).
The holomorphic $2$-differential $Q$ appearing in the $(2,0)$-part of $ds^2$ 
is defined on $M^2$, and is called the {\it Hopf differential\/} of $f$.
By definition, the umbilic points of $f$ equal the 
zeros of $Q$.  
Defining a meromorphic function on $\widetilde{M}^2$ by 
\begin{equation}\label{eq:rho-def}
   \rho=\theta/\omega \; , 
\end{equation}
then $|\rho|\colon{}M^2\to\R_+\cup\{0,\infty\}$ ($\R_+=\{r\in\R \,;\, r>0 \}$) is 
defined, and $p\in M^2$ is a singular point of $f$ if and only if
$|\rho(p)|=1$.

We note that the $(1,1)$-part of the first fundamental form 
\begin{equation}\label{eq:one-one-part}
     ds^2_{1,1}=|\omega|^2+|\theta|^2
\end{equation}
is positive definite on $M^2$, and 
$2ds^2_{1,1}$ coincides with the Sasakian metric's pull-back 
on the unit cotangent bundle $T^*_1 H^3$ by the Legendrian 
lift $L_f$ of $f$ (which is the sum of the first and third fundamental
forms in this case, 
see Section 2 of \cite{KUY2} for details).  
The two hyperbolic Gauss maps are 
\[
   G=\frac{A}{C}, \quad G_*=\frac{B}{D}, \quad
 \text{ where}\quad 
    \E_f=\begin{pmatrix} A & B \\ C & D \end{pmatrix}.
\]
Geometrically, $G$ and $G_*$ represent the intersection points in the 
ideal boundary $\partial H^3 = 
\C \cup \{ \infty \}$ of $H^3$ for the two oppositely-oriented normal geodesics 
emanating from $f$.  
The transformation $\E_f\mapsto a\E_f$ by 
$a=(a_{ij})_{i,j=1,2} \in \SL(2,\C)$ 
induces the rigid motion $f\mapsto afa^*$ as in \eqref{eq:isometric}, 
and $G$ and $G_*$ then change by the M\"obius transformation: 
\begin{equation}\label{eq:g-moebius}
    G\mapsto a\star G =  \frac{a_{11}G+a_{12}}{a_{21}G+a_{22}},\qquad
    G_*\mapsto a\star G_* =  \frac{a_{11}G_*+a_{12}}{a_{21}G_*+a_{22}}.
\end{equation}

\begin{remark}[The interchangeability of $\omega$ and $\theta$]
\label{rem:dual}
 The canonical forms $(\omega,\theta)$ have the $\U(1)$-ambiguity
 $(\omega,\theta)\mapsto (e^{is}\omega,e^{-is}\theta)$ ($s\in\R$)
 which corresponds to 
 \begin{equation}\label{eq:u1-amb}
    \E^{}_f\longmapsto 
    \E^{}_f
    \begin{pmatrix}
      e^{is/2} & 0 \\
      0  & e^{-is/2}
    \end{pmatrix}.
 \end{equation}
 For a second ambiguity, 
 defining the {\it dual} $\E_f^{\natural}$ of $\E_f$ by 
 \[
   \E_f^{\natural} = \E_f \begin{pmatrix}
                     0 & \iii \\ \iii & 0 
                    \end{pmatrix} \; , 
 \]
 then $\E_f^{\natural}$ is also Legendrian with 
 $f=\E_f^{\natural}{\E_f^{\natural}}^*$. 
 The hyperbolic Gauss maps $G^{\natural}$, $G_*^{\natural}$ 
 and canonical forms $\omega^{\natural}$, $\theta^{\natural}$ of 
 $\E_f^{\natural}$ satisfy 
 \begin{equation*}
  G^{\natural}=G_*, \quad G_*^{\natural} =G, \quad 
   \omega^{\natural}=\theta
   \quad\text{and}\quad
   \theta^{\natural}=\omega. 
 \end{equation*} 
 Namely, the operation $\natural$ interchanges the roles of $\omega$ and
 $\theta$ and also of $G$ and $G_*$.  
\end{remark}
The following fact holds (see \cite{KUY2} for fronts
and \cite{GMM} for the regular case):
\begin{fact}\label{fact:ode}
 Let $\omega$, $\theta$ be 
 holomorphic $1$-forms on a simply-connected Riemann surface $M^2$
 with $|\omega|^2+|\theta|^2$ positive definite.
 Solving the ordinary differential equation
 \[
    \E^{-1}d\E=
     \begin{pmatrix}
      0 & \theta \\
      \omega & 0
     \end{pmatrix},\qquad
   \E(z_0)=     \begin{pmatrix}
      1 & 0 \\
      0 & 1
     \end{pmatrix}
 \]
 gives a holomorphic Legendrian immersion of $M^2$ into $\SL(2,\C)$, 
 where $z_0\in M^2$ is a base point, and its projection into $H^3$ is 
 a flat front with canonical forms $(\omega,\theta)$.  
\end{fact}
\begin{remark}
 \label{rem:branch}
 If $|\omega|^2+|\theta|^2$ vanishes at a point
 $p\in M^2$, then $p$ is called a {\it branch point} of $f$.
 At such a branch point, $f$ is not a front, but the unit normal 
 vector field still extends 
 smoothly across $p$, so $f$ can be considered as a
 {\it frontal map}.
\end{remark}
\begin{remark}\label{rem:parallel}
 Considering $H^3$ as the hyperboloid in $\Lor^4$,
 the parallel front $f_t$ of $f$ is as in \eqref{eq:parallel}.
 As pointed out in \cite{GMM} and \cite{KUY2},
 \begin{equation}\label{eq:parallel-lift}
    \E^{}_{f_t}=
       \E^{}_f 
       \begin{pmatrix}
            e^{t/2} & 0 \\ 0 & e^{-t/2}
       \end{pmatrix}.
 \end{equation}
 Then the canonical forms $\omega_t$, $\theta_t$ 
 and the function $\rho_t=\theta_t/\omega_t$ of $f_t$ are written
 as
 \begin{equation}\label{eq:parallel-canonical}
    \omega_t = e^{t} \omega,\qquad
    \theta_t = e^{-t}  \theta\quad\text{and}\quad
    \rho_t = e^{-2t}\rho.
 \end{equation}
\end{remark}

\begin{fact}[\cite{KUY1}]
 \label{fact:legendre}
 For an arbitrary pair $(G,\omega)$ of a non-constant
 meromorphic function $G$ and a non-zero meromorphic $1$-form $\omega$ on $M^2$,
 the meromorphic map 
 \begin{equation}\label{eq:sm-analogue2}
  \E=\begin{pmatrix}
     A & dA/\omega \\
     C & dC/\omega 
    \end{pmatrix}\qquad
    \left(C=i\sqrt{\frac{\omega}{dG}},~
          A=GC\right)
 \end{equation}
 is a meromorphic Legendrian curve in $\PSL(2,\C)$ whose hyperbolic
 Gauss map and canonical form are $G$ and $\omega$,
 respectively.
 Conversely, if $\E$ is a meromorphic Legendrian curve in $\PSL(2,\C)$
 defined on $M^2$ with non-constant hyperbolic Gauss map $G$ and
 non-zero canonical form $\omega$, 
 then $\E$ is as in \eqref{eq:sm-analogue2}.
\end{fact}

\begin{remark}
 The $\E$ in \eqref{eq:sm-analogue2} has a sign 
 ambiguity, due to the square root of one meromorphic
 function.  So $\E$ is not defined in $\SL(2,\C)$, but 
 rather only in $\PSL(2,\C)$.  
\end{remark}

\begin{fact}[\cite{KUY1}]
 \label{fact:g-rep}
 Let $G$ and $G_*$ be non-constant meromorphic functions on a Riemann
 surface $M^2$ such that $G(p)\ne G_*(p)$ for all $p\in M^2$.
 Assume that 
\[
   \int_{\gamma} \frac{dG}{G-G_*} \in i \R 
\]
 for every loop $\gamma$ on $M^2$.  Set 
 \begin{equation}\label{eq:xi}
    \xi(z) = c \cdot \exp\int_{z_0}^{z}\frac{dG}{G-G_*},
 \end{equation}
 where $z_0\in M^2$ is a base point and $c\in \C\setminus\{0\}$ is
 an arbitrary constant.
 Then 
 \begin{equation}\label{eq:repG2}
     \E = \begin{pmatrix}
      G/\xi &      \xi G_{*}/(G-G_{*}) \\
      1/\xi &      \hphantom{G_*}\xi/(G-G_{*}) 
     \end{pmatrix}
 \end{equation}
 is a non-constant meromorphic Legendrian curve  
 defined on $\widetilde M^2$ in $\PSL(2,\C)$ whose hyperbolic Gauss maps
 are $G$ and $G_*$, and the projection $f=\E\E^*$  is
 single-valued on $M^2$. 
 Moreover, $f$ is a front if and only if $G$ and $G_*$ have no common
 branch points.  
 Conversely, any non-totally-umbilic flat front can be constructed 
 this way.
\end{fact}
\begin{remark}
 In \cite[Theorem 2.11]{KUY2}, we assumed that 
 all poles of $dG/(G-G_*)$ are of order $1$.
 This condition is satisfied automatically since $G(p)\neq G_*(p)$
 for all $p\in M^2$.
\end{remark}
If we write the constant $c$ in \eqref{eq:xi} as 
$c=e^{-(t+is)/2}$ ($t,s\in\R$), $s$ corresponds to the 
$\U(1)$-ambiguity \eqref{eq:u1-amb} and $t$ corresponds
to the parallel family \eqref{eq:parallel-lift}.  The canonical forms 
$\omega$, $\theta$ and Hopf differential $Q$ of $f$ in
Fact~\ref{fact:g-rep} are written as
\begin{equation} \label{eq:data-g}
  \omega = - \frac{1}{\xi^2} dG, \qquad 
  \theta = \frac{\xi^2}{(G-G_*)^2} dG_*,\qquad
 Q = \frac{-dG\,dG_*}{(G-G_*)^2}.
\end{equation}

Let $z$ be a local complex coordinate on $M^2$ and write
$\omega=\hat\omega\,dz$ and $\theta=\hat\theta\,dz$.
Then we have the following identities (see \cite{KUY2}):
\begin{align}
  &\frac{\hat\omega'}{\hat\omega}=\frac{G''}{G'}-2\frac{G'}{G-G_*}, \quad
  \frac{\hat\theta'}{\hat\theta}=\frac{G''_*}{G'_*}-2\frac{G'_*}{G_*-G},
  \label{eq:can-der-Gauss}
  \\
  &s(\omega)=2Q+S(G), 
  \quad s(\theta)=2Q+S(G_*),
  \label{eq:schwarz}
\end{align}
 where ${}'=d/dz$, and $S(G)$ is 
the Schwarzian derivative of $G$ with respect to $z$ as in 
\begin{equation}\label{eq:schwarz2}
    S(G) = 
    \left\{\left(
                     \frac{G''}{G'}
                    \right)'-
    \frac{1}{2}
    \left(
     \frac{G''}{G'}
    \right)^2\right\}dz^2,
\end{equation}
 and $s(\omega)$ is the Schwarzian derivative of
 the integral of $\omega$, that is,
\begin{equation}\label{eq:small-s}
 s(\omega) = S(\varphi)
  =\left\{
 \left(\frac{\hat\omega'}{\hat\omega}\right)' 
 -\frac{1}{2}\left(\frac{\hat\omega'}{\hat\omega}\right)^2\right\}dz^2 
 \quad
 \left(\varphi(z)=\int_{z_0}^{z}\omega  \right).
\end{equation}
Note that although the Schwarzian derivative depends on the choice 
of local coordinates, 
the difference $S(G)-S(G_*)$ does not. 
If $G$ expands as 
$G(z)=a+b (z-p)^m+o\bigl((z-p)^m\bigr)$, $b\ne 0$, 
$m \in \Z_+$,
where $o\bigl((z-p)^m\bigr)$ denotes higher order terms, then 
\begin{equation}\label{eq:s_ord}
 S(G)=
  \frac{1}{(z-p)^2}\left(
			\frac{1-m^2}{2}+o(1)
		       \right)\,dz^2.
\end{equation}
Similarly, if a meromorphic $1$-form $\omega=\hat\omega\,dz$
has an expansion 
\begin{equation}\label{expandomega}
    \hat\omega(z) = c(z-p)^{\mu}\bigl(1+o(1)\bigr) 
  \qquad (c\neq 0,\mu\in\R),
\end{equation}
then we have
\begin{equation}\label{eq:s-order}
 s(\omega)
   =
   \frac{1}{(z-p)^2}\left(
       -\frac{\mu(\mu+2)}{2}
            +o(1)\right)\,dz^2.
\end{equation}
Conversely, if $\omega$ satisfies \eqref{eq:s-order}, it
expands as in \eqref{expandomega}.
For later use, we define the order of a metric defined on
a punctured disc.  

\begin{definition}[\cite{Troyanov}] \label{def:order}
 A conformal metric $d\sigma^2$ on the punctured disc 
 $D^*=\{z\in\C\,;\,0<|z|<1\}$
 is of {\em finite order\/} if there exist a $c>0$ and 
 $\mu\in \R$ so that $d\sigma^2$ is locally expressed as 
 \[
    d\sigma^2=c|z|^{2\mu} \bigl(1+o(1)\bigr) \, |dz|^2 \; . 
 \]
 We define the {\em order\/} $\ord_0 d\sigma^2$ of 
 $d\sigma^2$ at the origin to be $\mu$.  
\end{definition}

Finally, 
we give the following proposition on our definition of flat fronts:
\begin{proposition}
 Let $f:M^2\to H^3$ be a front such that the regular set
 is open and dense in $M^2$.
 Then the following two conditions are equivalent.
\begin{enumerate}
 \item\label{item:add:1} 
       For each point $p\in M^2$, there exists a
       real number $t_0$ such that $p$ is a 
       regular point of the parallel surface $f_{t_0}$
       and the Gaussian curvature of $f_{t_0}$ vanishes.
 \item\label{item:add:2} 
       The Gaussian curvature vanishes
       on the regular set of $f$.
\end{enumerate}
\end{proposition}
The first condition is in fact the definition of flat fronts.
\begin{proof}
 If $f$ satisfies \ref{item:add:1}, then  
 $f_t$ is also a flat front for each $t\in \R$. (See \cite{KUY2}.)
 Thus \ref{item:add:1} implies \ref{item:add:2} obviously.
 Next, we suppose $f$
 satisfies \ref{item:add:2}.
 It suffices to discuss under the assumption that $p$ is a singular 
 point of $f$.
 Note that $f$ has a smooth unit normal vector field $\nu$
 like as in the case of an immersion. 
 Now we shall show that
 the parallel front $f_t$ is an immersion at $p$
 for sufficiently small $t$.
 First, we consider the case that $df$ vanishes at
 $p$. Since $f$ is a front, $\nu:U\to \Lor^4$ is an immersion
 and so $f_t$ gives an immersion at $p$ for all $t\ne 0$.
 Next we consider the case that the kernel of $df$ at $p$ 
 is one dimensional, and take a non-vanishing vector $\eta\in T_pM^2$ 
 such that $df(\eta)=0$. 
 Then we can take a local 
 coordinate system $(U;u,v)$ centered at $p$ such that
 $\partial/\partial u=\eta$.
 Since $f_u:=df(\partial/\partial u)=0$, we have
 \begin{equation}\label{eqn:add}
  0=\inner{f_u}{\nu_v}=-\inner{f_{uv}}{\nu}=
   \inner{f_v}{\nu_u},
 \end{equation}
 where  $\inner{\,\,}{\,\,}$ is the canonical Lorentzian metric on $\Lor^4$.
 Since $f$ is a front and $f_u=0$, we have $f_v\ne 0$ and $\nu_u\ne 0$.
 Then \eqref{eqn:add} implies that $f_v$ and $\nu_u$ are linearly
 independent. Then by \eqref{eq:parallel}, 
 $f_t$ is an immersion at $p$ for sufficiently small $t$.
 (Moreover,  
 $f_t$ is an immersion for $t\in \R$ except for only one value.
 See Remark~\ref{rmk:add} below.)

 Let $K_t$ be the Gaussian curvature of $f_t$ ($t\ne 0$)
 near $p$. 
 Suppose that $K_{t_0}(p)\ne 0$ for a sufficiently small 
 $t_0$. Then any point $q$ near $p$ satisfies $K_{t_0}(q)\ne 0$.
 On the other hand, since the regular set of $f$ is dense in $M^2$,
 there exists a point $q$ sufficiently near $p$ such that
 the Gaussian curvature of $f$ vanishes 
 around $q$. Then $K_{t}$ vanishes as long as $q$ is a regular point
 of $f_t$, and we have $K_{t_0}(q)=0$, a contradiction. 
 Thus we have $K_{t_0}(p)=0$, which implies \ref{item:add:1}.
\end{proof}

\begin{remark}\label{rmk:add}
 Let $f:M^2\to H^3$ be a front and $p\in M^2$ 
 a regular point.
 Let $\lambda_j$ ($j=1,2$) be the principal curvatures of $f$ at $p$.
 Then, the parallel front $f_t$ has a singularity at $p$
 if and only if $\lambda_1=\coth t$
 or $\lambda_2=\coth t$.
 Suppose now that $f$ is flat and $p$ is a non-umbilical point. 
 Since $\lambda_1\lambda_2=1$, we may assume that 
 $\lambda_1<1<\lambda_2$.
 Then $p$ is a singular point of $f_t$ only when 
 $\lambda_2=\coth t$, and such a $t$ is uniquely determined.
\end{remark}

\section{Completeness and weak completeness}
\label{sec:complete}

Let $f\colon M^2 \to H^3$ be a flat front. We say that $f$ is 
{\it complete\/} if there exists a symmetric $2$-tensor field $T$ with
compact support
so that the sum $ T+ ds^2 $ is a complete Riemannian metric on 
$M^2$  (cf.~\cite{KUY1}).  (Because this metric is required to be Riemannian, if the 
singular set accumulates at some end of $f$, then, by definition, $f$ is 
not complete.)  We say that $f$ is {\it weakly complete\/} if 
the $(1,1)$-part $ds^2_{1,1} = |\omega|^2 +|\theta|^2$
in \eqref{eq:one-one-part} of the induced metric
is complete and Riemannian on $M^2$. 
Since $ds^2_{1,1}$ is proportional to the pullback of the Sasakian metric, 
this definition of weak completeness is analogous to the notion of completeness
in Melko and Sterling \cite{MS} for constant negative curvature surfaces in $\R^3$.

We say that a flat front $f$ is of {\em finite type\/} 
if $ds^2_{1,1}$ has finite total curvature.  
A 2-manifold $M^2$ is said to have {\it finite topology\/} if 
there exist a compact $2$-manifold $\overline{M}^2$ and finitely many points 
$p_1, \dots , p_n \in \overline{M}^2$ such that 
$M^2$ is homeomorphic to $\overline{M}^2 \setminus \{p_1 , \dots , p_n\}$.  
A small neighborhood $U_j$ of a puncture 
point $p_j$, or even just the puncture point $p_j$ itself, 
is called an {\it end} of $f$.  
An end $p_j$ is {\em complete\/}, resp. {\em weakly complete\/}, 
if $ds^2$, resp. $ds^2_{1,1}$, is complete at $p_j$.
\begin{proposition}\label{prop:f->finite}
 A complete flat front is weakly complete and of finite type. 
\end{proposition}
\begin{proof}
 Let $f\colon M^2 \to H^3$ be a complete flat front, then 
 $f$ is weakly complete by \cite[Corollary 3.4]{KUY2}.
 We now show that $ds^2_{1,1}$ for $f$ has finite total curvature: 
 By \cite[Lemma 3.3]{KUY2}, $M^2$ is biholomorphic to 
 $\overline{M}^2 \setminus \{p_1, \dots , p_n\}$ for some 
 compact Riemann surface $\overline{M}^2$. 
 By completeness, there exist 
 neighborhoods $U_j$ of each $p_j$ so that $U_j \setminus \{p_j \}$ 
 contains no singularities.  
 Hence \cite[Lemma 2]{GMM} implies 
 \begin{equation}\label{canonicalformforomegandtheta}
   \omega = z^{\mu} \omega_0, \quad 
   \theta = z^{\nu} \theta_0 \quad 
   (\mu, \nu \in \R), 
 \end{equation}
 where $z$ is a local coordinate with $z(p_j)=0$, 
 and $\omega_0, \theta_0$ are single-valued 
 holomorphic $1$-forms on $U_j$ which are nonzero at $p_j$.  
 Thus, the order of $ds^2_{1,1}$ at $p_j$ is finite 
 (see Definition~\ref{def:order}).
 Recalling the formula 
 (cf.\ \cite{Fang}, \cite{Shiohama}) 
 \begin{equation}\label{eq:total-curvature}
    \frac{1}{2\pi} \int_{M^2} (-K_{ds^2_{1,1}}) dA 
     = -\chi(\overline{M}^2) 
       + \sum_{j=1}^n \ord_{p_j} ds^2_{1,1},
 \end{equation}
 where $K_{ds^2_{1,1}}$ and $dA$ denote the Gaussian curvature and
 area element of $(M^2, ds^2_{1,1})$, and 
 $\chi(\overline M^2)$ the Euler number of $\overline{M}^2$, 
 we see that $ds^2_{1,1}$ has finite total curvature.  
\end{proof}

\begin{proposition}\label{prop:finite->}
 When $f\colon M^2 \to H^3$ is a weakly complete flat front 
 of finite type, 
\begin{enumerate}
 \item\label{item:finite-1}
      $M^2$ is biholomorphic to 
      $\overline{M}^2 \setminus \{p_1, \dots , p_n \}$,  
      for some compact Riemann surface $\overline{M}^2$ and finitely 
      many points $p_1, \dots, p_n \in \overline M^2$, 
 \item\label{item:finite-2} 
      $ds^2_{1,1}$ has finite order at each $p_j$, and the 
      canonical $1$-forms 
      $\omega$, $\theta$ are of finite order, and 
 \item\label{item:finite-3} 
      $Q$ is a meromorphic differential on $\overline{M}^2$.  
\end{enumerate}
\end{proposition}
\begin{proof}
\ref{item:finite-1}:
   By Huber's theorem,  
   $M^2$ is diffeomorphic to 
   $\overline{M}^2 \setminus \{p_1, \dots , p_n \}$ 
   since $ds^2_{1,1}$ is complete with finite total curvature. 
   In fact, they can be biholomorphic, as the Gaussian curvature 
   of $ds^2_{1,1}$ satisfies $K_{ds^2_{1,1}} \leq 0$ (\eqref{Krelations} below 
   implies $K_{ds^2_{1,1}} \leq 0$).  
\par
\ref{item:finite-2}:
   We shall show that each of $|\omega|^2$, $|\theta|^2$ has 
   finite order at $p_j$. 
   Take a local coordinate $z$ such that $z(p_j)=0$. 
   Since $|\omega|$, $|\theta|$ are both single-valued on $M^2$, 
   there exist real numbers $\mu, \nu \in [0,1)$ such that 
   $\omega \circ \tau = e^{2 \pi \iii \mu} \omega$ and 
   $\theta \circ \tau = e^{2 \pi \iii \nu} \theta$ 
   for the deck transformation $\tau$ 
   associated to a loop wrapped once about $p_j$. 
   Thus 
 \begin{equation}\label{eq:omega-0}
    \omega = z^{\mu} \omega_0 \quad\text{and}\quad
    \theta = z^{\nu} \theta_0,  
 \end{equation}
  where $\omega_0, \theta_0$ are single-valued 
  holomorphic $1$-forms on a punctured neighborhood
  $D^*(\varepsilon)=\{ 0 < |z| < \varepsilon\}$. 
  The function $\rho$ in \eqref{eq:rho-def} can be written as 
 \begin{equation}\label{eq:rho-0}
    \rho = \frac{\theta}{\omega} = z^{\mu-\nu}\frac{\theta_0}{\omega_0}
         = z^{\mu-\nu} \rho_0 \; ,  
 \end{equation}
  where $\rho_0=\theta_0/\omega_0$ is a single-valued holomorphic function on
  $D^*(\varepsilon)$. 

 First, we show that $\rho_0$ has at most a pole at $z=0$, that is, not an 
 essential singularity.  
 Consider a constant mean curvature $1$ surface 
  \[
    f_1\colon{}\widetilde{D^*(\varepsilon)}\longrightarrow H^3
 \]
 of the universal cover $\widetilde{D^*(\varepsilon)}$ into $H^3$ 
 with Weierstrass data 
 $(g_1,\omega_1)=(\rho,\hat\omega^2\,dz)$, where
 $\omega=\hat\omega\,dz$ (see \cite{UY}).
 Since the induced metric $ds^2_1$ by $f_1$ and $ds^2_{1,1}$ are 
 \[
     ds^2_1 = h^2\,|dz|^2,\qquad
     ds^2_{1,1}= h\,|dz|^2,\quad h=(1+|\rho|^2)|\hat\omega|^2 ,
 \]
$ds_{1,1}^2$ is positive definite, so $f_1$ is an immersion.  Also, 
 \begin{equation}\label{Krelations}
    K_{ds^2_1} dA_{ds^2_1} = 2 K_{ds^2_{1,1}} dA_{ds^2_{1,1}},
\end{equation}
 holds, where 
 $K_{ds^2_1}$ (resp.~$K_{ds^2_{1,1}}$) and $dA_{ds^2_1}$ (resp.~$dA_{ds^2_{1,1}}$) 
 are the Gaussian curvature and area element 
 with respect to the metric $ds^2_1$ (resp.~$ds^2_{1,1}$).  
 The induced metric $ds^2_{1}$ is well-defined on
 $D^*(\varepsilon)$, because $|\omega|$ and $|\theta|$ are single-valued.
 Since $f$ is of finite type, the total curvature
 \[
    \int_{D^*(\varepsilon)} K_{ds^2_1} dA_{ds^2_1}
 \]
 is finite.  Then \cite[Proposition 4]{Br} implies 
 $\rho_0$ in \eqref{eq:rho-0} has at most a pole at $z=0$.  So 
 if $\omega_0$ in 
 \eqref{eq:omega-0} has at most a pole at $z=0$, the same is true of $\theta_0$ 
 and (2) will be proven.  
 Taking the dual as in Remark~\ref{rem:dual} if need be, 
 we may assume $|\rho(0)|<\infty$.
 In a sufficiently small neighborhood of $z=0$, we have
 \[
    ds^2_{1,1} = (1+|\rho|^2)|\omega|^2 
     \leq k_1 |\omega|^2 
      = k_1 |z|^{2\mu}|\omega_0|^2 
         \leq k_2 |\omega_0|^2 
 \]
 for some constants $k_1,k_2>0$, since $\mu\in [0,1)$.
 Completeness of $ds_{1,1}^2$ implies 
 $k_2 |\omega_0|^2$ is also complete at $z=0$. 
 Hence, $\omega_0$ has at most a pole at $z=0$ (\cite[Lemma 9.6]{O}).  
\par
\ref{item:finite-3}:
 Since $Q=\omega\theta$, assertion (3) is immediate from 
 the proof of \ref{item:finite-2}. 
\end{proof}
Propositions~\ref{prop:f->finite} and  \ref{prop:finite->}
yield:   
\begin{theorem}\label{thm:f<->finite-cs}
 A flat front $f\colon M^2 \to H^3$ is complete if and only if 
 \begin{enumerate}
  \item\label{item:finite-thm-1} 
       $f$ is weakly complete and of finite type, and 
  \item\label{item:finite-thm-2} 
       the set of singularities, denoted by 
       $\Sigma_f \subset M^2$,  is compact. 
 \end{enumerate}
\end{theorem}
\begin{proof}
 Proposition \ref{prop:f->finite} shows that completeness implies 
 \ref{item:finite-thm-1}, \ref{item:finite-thm-2}.  
 Now suppose \ref{item:finite-thm-1}, \ref{item:finite-thm-2} hold. 
 Proposition \ref{prop:finite->} implies $M^2$ is biholomorphic to 
 $\overline M^2 \setminus \{ p_1 ,\dots , p_n\}$, and that $f$'s 
 canonical $1$-forms $\omega$, $\theta$ are of finite order.  
 Since $|\omega|$ and $|\theta|$ 
 are well-defined, $\omega$ and $\theta$ are in the form 
 \eqref{canonicalformforomegandtheta}, 
 and at least one of $|\omega|^2$, $|\theta|^2$ is complete, at 
 any $p_j$.  If $|\rho(p_j)| = 1$ for $\rho=\theta/\omega$ at some $p_j$, 
 then $\rho$ would be locally holomorphic at $p_j$ and (2) would not hold, 
 so $|\rho(p_j)| \ne 1$ at every $p_j$.  
 To prove completeness of $f$, 
 we must show $ds^2 = |\omega + \bar \theta|^2$ is complete 
 at all $p_j$. 
 Without loss of generality, we may assume $\ord_{p_j}|\theta|^2
 \ge \ord_{p_j}|\omega|^2$, so $|\omega|^2$ is complete at $p_j$.  
 So 
 \[
   ds^2=|\omega+\bar \theta|^2\ge 
   \bigl||\omega|-|\theta|\bigr|^2=|\omega|^2 
   \bigl|1-|\rho|\bigr|^2.
 \]
 Since $|\rho(p_j)| \not\in \{ 1,\infty \}$, it follows that 
 $ds^2$ is complete at $p_j$.  
\end{proof}
Proposition \ref{prop:f->finite} and the following theorem
prove the introduction's Theorem~\ref{thm:A}.
\begin{theorem}\label{thm:fff->f_t}
 Let $f$ be a weakly complete flat front of finite type, 
 with $n$ ends.  Then 
 the parallel fronts $f_t$ of $f$ are complete 
 except for at most $n$ values of $t$.
\end{theorem}
\begin{proof}
 By Theorem~\ref{thm:f<->finite-cs}, 
 we need only show compactness of the set 
 of singularities $\Sigma(t)$ of $f_t$ away from 
 at most $n$ values of $t$. 
 Since $f$ is weakly complete and of finite type, 
 we can set $M^2=\overline{M}^2\setminus\{p_1,\dots,p_n\}$ 
 for some compact Riemann surface $\overline{M}^2$, 
 and $\omega$, $\theta$ are both of finite order at each $p_j$, 
 by Proposition~\ref{prop:finite->}.  Then the function 
\[
   |\rho|=\left|\frac{\theta}{\omega}\right|:
   \overline{M}^2\longrightarrow \R_+\cup\{0,\infty\}
\]
 is well-defined and continuous.  The closure of the singular set 
 $\overline{\Sigma(t)}$ in $\overline{M}^2$ is 
 \[
    \overline{\Sigma(t)}=\{p\in \overline{M}^2\,;\,
         |\rho(p)|=e^{2t}
         \}.
 \]
 Thus $\Sigma(t)$ is compact when 
 $\{p_1,\dots, p_n\}\cap \overline{\Sigma(t)}$ 
 is empty.  
 Let $\{p_{j_1},\dots,p_{j_m}\}$ be the subset of ends such that
 $|\rho(p_{j_k})|\neq 0,\infty$.
 Taking the unique $t_{j_k}\in\R$
 so that $p_{j_k}\in \overline{\Sigma(t_{j_k})}$ for each $k$,
 i.e. $|\rho(p_{j_l})|=\exp(2t_{j_k})$, 
 $\Sigma(t)$ is compact for any 
 $t\in \R \setminus \{t_{j_1},\dots,t_{j_m} \} $.
\end{proof}
\begin{figure}
   \begin{tabular}{c@{\hspace{1cm}}c@{\hspace{1cm}}c}
        \includegraphics[width=3cm]{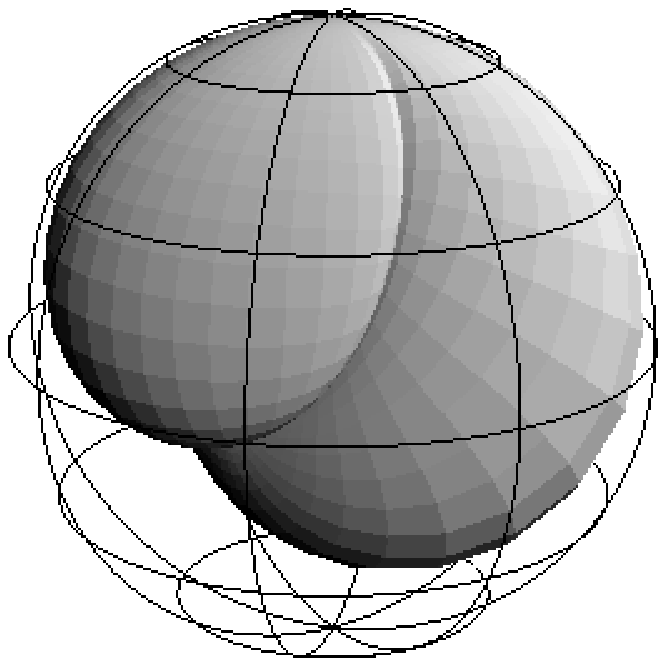} &
        \includegraphics[width=3cm]{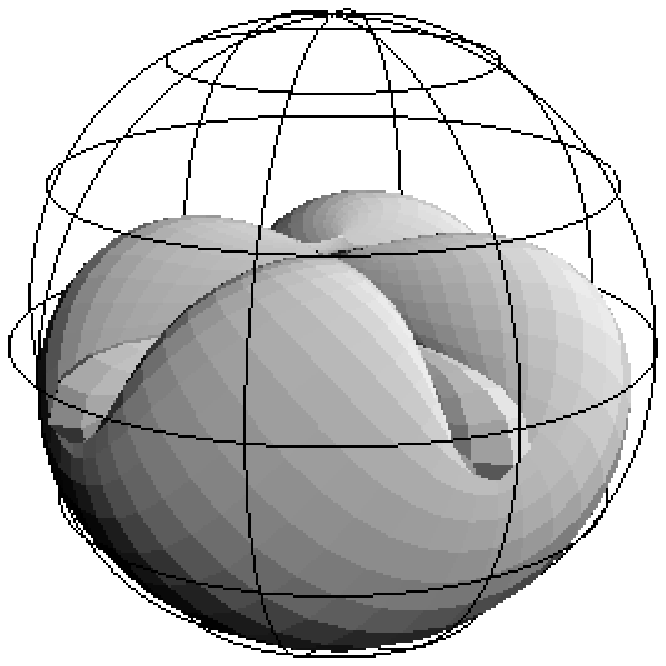} &
        \includegraphics[width=3cm]{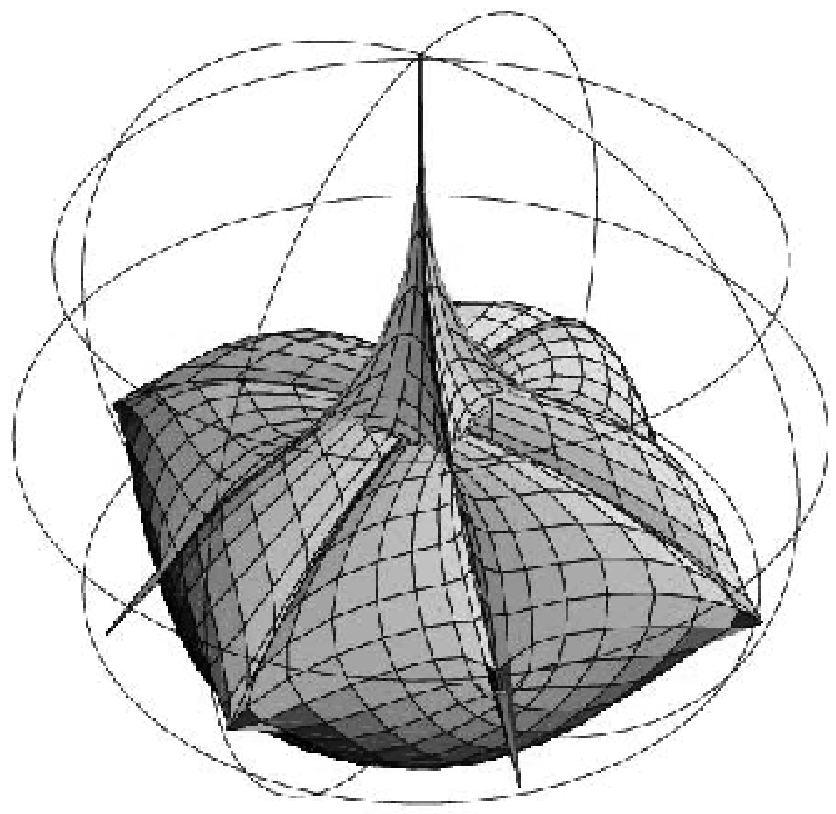} \\
  \end{tabular}
\caption{%
 A peach front (Example~\ref{ex:peach}) on the left, and a genus 1 
 complete front with 5 embedded ends in the middle (Example~\ref{ex:4k+1}) 
 with its caustic on the right.
}
\label{fig:peach}
\end{figure}
\begin{definition}\label{def:multiplicity}
 Let $f\colon{}\overline{M}^2\setminus \{p_1,\dots,p_n\}\to H^3$
 be a weakly complete flat front, 
 and $(U,z)$ a complex coordinate of $\overline{M}^2$ 
 with $z(p_j)=0$.  
 Suppose that all ends are {\it regular\/}, i.e. 
 $G$ and $G_*$ have at most poles at 
 $p_1,\dots,p_n$, and that both $G$ and $G_*$ are non-constant.  
 By a suitable motion in $H^3$, we may assume 
 $G$ and $G_*$ have no poles at $p_j$.  Then we have, with $m,m_* \in \Z_+$, 
 \begin{equation}\label{eq:g-expansion}
  G(z)=a+b z^m+o(z^m)
    \quad\text{and}\quad
  G_*(z)=a_*+b_* z^{m_*}+o(z^{m_*}),
 \end{equation}
 where $a,a_*\in\C$, $b,b_*\in \C\setminus\{0\}$
 and $o(z^m)$ and $o(z^{m_*})$ are 
 higher order terms.  We set
 \[
    r_{p_j}(G)=m,\qquad  r_{p_j}(G_*)=m_*
 \]
 to be the {\em ramification numbers\/} of $G$ and $G_*$ at $p_j$, 
 respectively.  Moreover, we define the {\it multiplicity\/} of $f$ at the 
 end $p_j$ to be 
 \[
     m(f,p_j)=\min\{m,m_*\}
     \quad\bigl(=\min\{r_{p_j}(G),r_{p_j}(G_*)\}\bigr) \; . 
 \]
\end{definition}

If $f$ is complete at an end $p$, 
{\it $m(f,p)=1$ if and only if the end $p$ is properly embedded} \cite{KUY2}.  
Roughly speaking, the multiplicity 
of a complete end is the winding number of a slice of the end
(see \cite{KRUY}).  We have the following:
\begin{proposition}\label{prop:ordQgem1}
 Let $f\colon{}\overline{M}^2\setminus \{p_1,\dots,p_n\}\to H^3$
 be a weakly complete flat front whose ends are all regular.
 If $Q$ has at most a simple pole at an end $p_j$, 
 then $|\omega|^2$ and $|\theta|^2$ have finite orders at $p_j$ and 
 the following identity holds: 
 \begin{equation}\label{ram-number-multiplicity}
  m(f,p_j)=
   \min\left\{\bigl|1+\ord_{p_j}|\omega|^2\bigr|,
              \bigl|1+\ord_{p_j}|\theta|^2\bigr|
             \right\} \; . 
 \end{equation}
\end{proposition}
\begin{proof}
 Using a complex coordinate $z$ with $z(p_j)=0$, 
 assume $G$ and $G_*$ expand as in \eqref{eq:g-expansion}.  As 
 $Q$ has at most a simple pole at $p_j$, 
 \eqref{eq:schwarz} and \eqref{eq:s_ord} give 
 \[
    s(\omega) = 
      \frac{1}{z^2}
      \left(
        \frac{1-m^2}{2}+o(1)
      \right)\,dz^2,\qquad
    s(\theta) = 
      \frac{1}{z^2}
      \left(
        \frac{1-m_*^2}{2}+o(1)
      \right)\,dz^2.
 \]
 So, by Section \ref{sec:prelim}, 
 $|\omega|^2$ and $|\theta|^2$ are of finite order (and well defined) 
 at $p_j$.  Thus $\omega$ and $\theta$ are of the form in 
 \eqref{canonicalformforomegandtheta}, and 
 \[
    \mu=\ord_{p_j}|\omega|^2=\pm m-1,\qquad
    \mu_*=\ord_{p_j}|\theta|^2=\pm m_*-1.
 \]
 Hence
 $m(f,p_j) = \min\{m,m_*\}=\min\left\{|\mu+1|,|\mu_*+1| \right\}$ satisfies 
 \eqref{ram-number-multiplicity}.  
\end{proof}

Next, we give a weakly complete example that is neither 
complete nor of finite type.

\begin{example}[The peach front in \cite{KUY2}] \label{ex:peach}
 Let $b\in \C$ be a non-vanishing complex number.
 We define rational functions on $\C\cup\{\infty\}$ by
 \[
   G=z,\qquad  G_*=z-b.
 \]
 By Fact~\ref{fact:g-rep}, we get a holomorphic Legendrian curve
 \[
     \E=\left(\begin{array}{rr}
          \dfrac{z}{c}\,e^{-z/b} & \dfrac{c}{b} (z-b) e^{z/b} \\[6pt]
          \dfrac{1}{c}\,e^{-z/b} & \dfrac{c}{b} e^{z/b}
        \end{array}\right)
 \]
 with 
 \[
    \omega=-\frac{1}{c^2}e^{-2z/b}\,dz, \qquad 
    \theta=\frac{c^2}{b^2}e^{2z/b}\,dz, \qquad 
   ds^2_{1,1}=|\omega|^2+|\theta|^2\ge \frac{2}{|b|^2}|dz|^2,
 \]
 which implies that $f=\E\E^*$ is a weakly complete
 flat front in $H^3$, called a peach front. 
 As the singular set of $f$ given by 
 $|\theta|=|\omega|$ accumulates at $z=\infty$ for all $c$, $f$ and its 
 parallel fronts are all not complete.
 By Theorem~\ref{thm:fff->f_t}, $f$ is not of finite type.
 As we will see in Section~\ref{sec:caustic}, the peach fronts are 
 characterized by the property that their caustics are horospheres.  
 Figure~\ref{fig:peach} shows the peach front for $b=c=1$.
\end{example}

\section{Examples}
\label{sec:example}
Let $M^2$ be a Riemann surface and $d\sigma^2$ 
a flat metric compatible to the conformal structure of $M^2$
(we call $d\sigma^2$ a ``flat conformal metric'').
Then there exists a developing map (as a holomorphic map)
\[
   \varphi\colon{}\widetilde{M}^2\longrightarrow \R^2=\C
\]
such that $d\sigma^2$ is the pull-back of the canonical
metric of $\R^2$, where $\pi\colon{}\widetilde{M}^2\to M^2$ is the
universal cover of $M^2$.  The differential $d\varphi$ is 
a holomorphic $1$-form on $\widetilde M^2$, and 
\begin{equation}\label{eq:developing}
  d\sigma^2=|d\varphi|^2.
\end{equation}
Such a $1$-form $d\varphi$ is determined up to multiplication by 
a unit complex number.  
We call $d\varphi$ the {\em associated $1$-form\/} of the metric
$d\sigma^2$.

For example, a flat front $f\colon{}M^2\to H^3$ 
without umbilics gives 
two flat conformal  metrics $|\omega|^2$ and $|\theta|^2$ 
globally defined on $M^2$, with associated
 $1$-forms $\omega$ and $\theta$.  
(At an umbilic point $q$ of $f$, one of $\omega$ and $\theta$ will vanish, 
see \eqref{fff-sff}. 
So one of the metrics $|\omega|^2$ or $|\theta|^2$ will degenerate at $q$.)

To construct examples, we introduce the following result.
\begin{theorem}\label{thm:met}
 Let $M^2$ be a Riemann surface and $|\omega|^2$ be a 
 flat conformal metric on $M^2$ with associated $1$-form $\omega$.
 Let $G$ be a meromorphic function on $M^2$.  
 Suppose that $d\sigma^2 = |\omega|^2+|\theta|^2$ is a smooth and positive
 definite metric on $M^2$, where 
 \begin{equation}\label{eq:theta}
    \theta=\frac{Q}{\omega},\qquad 
     Q=\frac{s(\omega)-S(G)}{2}.
 \end{equation}
 Then the map $f=\E\E^*:M^2\to H^3$ given by \eqref{eq:sm-analogue2}
 gives a flat front with 
 canonical forms $\omega$ and $\theta$,
 Hopf differential $Q$, and $G$ as one of its hyperbolic Gauss maps. 
\end{theorem}
Theorem~\ref{thm:met} follows directly from Fact~\ref{fact:legendre}
and \eqref{eq:schwarz}. 
Moreover, we have: 
\begin{proposition}\label{prop:order}
 In the situation of Theorem~\ref{thm:met}, suppose 
 also that $M^2$ is biholomorphic to 
 $\overline{M}^2\setminus \{p_1,\dots,p_n\}$, 
 for some compact Riemann surface $\overline{M}^2$.  
 Then each end $p_j$ of $f$ is 
 weakly complete if $p_j$ is a pole of $Q$ of order $2$.
\end{proposition}
\begin{proof}
If $Q=\omega\theta$ has a pole of order $2$ at $p_j$, we have
\[
   \ord_{p_j}|\omega|^2+\ord_{p_j}|\theta|^2=-2.
\]
Thus $\min\{\ord_{p_j}|\omega|^2,\ord_{p_j}|\theta|^2\}\le -1$
and $ds^2_{1,1}$ is complete.
\end{proof}

Here we construct 
higher genus flat fronts with regular ends, 
for which the ends might not be embedded.  
We begin with the following result, due to Troyanov:  

\begin{fact}\label{fact:Tro}
 Let $\overline{M}^2$ be a compact Riemann surface with 
 $p_1,\dots,p_n\in {\overline M}^2$,
 and let $\mu_1,\dots,\mu_n$  be real numbers
 which satisfy
 \[
    \chi({\overline M}^2)+\sum_{j=1}^n \mu_j=0.
 \]
 Then there exists a flat conformal {\rm (}positive definite{\rm )}
 metric $d\sigma^2$ 
 on $M^2=\overline{M}^2\setminus\{p_1,\dots,p_n\}$ 
 such that $\ord_{p_j}d\sigma^2=\mu_j$ 
 for each $j=1,\dots,n$. 
 Such a metric is unique up to homothety.  
\end{fact}
This fact is proved in \cite{Troyanov} for $\mu_j>-1$, 
but the same argument works for any real numbers $\mu_j$.
For the metric $d\sigma^2$ in Fact~\ref{fact:Tro}, the formal sum
$\mu_1 p_1+\dots+\mu_np_n$ 
is called the {\em singular divisor\/} of $d\sigma^2$.  
With it, we can prove the following: 
\begin{theorem}\label{thm:troyanov-example}
 Given any compact Riemann surface $\overline{M}^2$ and
 meromorphic function $G$ on $\overline{M}^2$ with 
 branch points $p_1,\dots,p_n$, there exists a complete flat front 
 $f\colon{}\overline{M}^2 \setminus \{p_1,\dots,p_n\}\to H^3$
 with all ends $p_j$ regular and with hyperbolic Gauss map $G$.  
\end{theorem}

\begin{proof}
 By a suitable motion in $H^3$, we may assume $G$
 has no poles at $p_1,\dots,p_n$.
 Let $m_j$ be the ramification number of $G$ at $p_j$
 (see Definition~\ref{def:multiplicity}).
 We can choose an $n$-tuple of real numbers
 $\mu_1, \dots , \mu_n$ 
 such that $\mu_j \neq \pm  m_j-1$  $(j=1,\dots,n)$ and
 \[
   \chi({\overline M}^2)+\sum_{j=1}^n \mu_j=0 \; . 
 \]
 By Fact~\ref{fact:Tro}, there exists a flat conformal metric
 $d\sigma^2$ on $M^2=\overline{M}^2\setminus\{p_1,\dots,p_n\}$
 with singular divisor $\mu_1p_1+\cdots +\mu_n p_n$.
 Then we can write 
 $d\sigma^2=|\omega|^2$ where $\omega$ is a holomorphic 
 $1$-form on $\widetilde{M}^2$, see \eqref{eq:developing}.
 By Theorem~\ref{thm:met} we can construct a flat front 
 $f\colon{}\widetilde{M}^2\to H^3$ from the data $(G,\omega)$. 
 Moreover, since $|\omega|$ is well-defined on $M^2$, 
 any given deck transformation of $M^2$ changes $\omega$ to 
 $e^{i \beta} \omega$ for some $\beta \in \R$.  Hence 
 $f$ is single-valued on $M^2$ by 
 \eqref{eq:sm-analogue2}, that is, $f\colon{}M^2\to H^3$.
 Since $\mu_j \neq \pm  m_j-1$, $Q$ has a pole of order 
 $2$ at each $p_j$, by \eqref{eq:schwarz}.  
 Then by Proposition~\ref{prop:order}, 
 $f$ is weakly complete at $p_j$.  
 Then, since $\ord_{p_j}|\omega|^2$ and $\ord_{p_j}Q$ are finite, 
 $\ord_{p_j}|\theta|^2$ is as well, so $f$ is of finite type.  
 Then by Theorem~\ref{thm:fff->f_t},
 there exist infinitely many complete flat fronts in the  parallel
 family
 of $f$, all with hyperbolic Gauss map $G$.  
 Since $G$ has no essential singularity, each $p_j$ is 
 a regular end (see Proposition \ref{prop:regularend} below).  
\end{proof}

\begin{remark}
 Let $\overline{M}^2$ be a genus $k$ hyperelliptic Riemann surface 
 with associated meromorphic function 
 $G\colon{}\overline{M}^2 \to \C\cup \{ \infty \}$ of degree $2$ 
 having branch points $p_1,\dots,p_n$ ($n=2k+2$).
 Take a suitable integer $m$ and reals $\nu_1,\dots,\nu_m<-1$ so that 
 \[
      \chi(\overline{M}^2)+n+\sum_{j=1}^m \nu_j=0.
 \]
 Choosing 
 points $q_1,\dots,q_m\in \overline{M}^2\setminus\{p_1,\dots,p_n\}$,
 there exists a flat conformal metric 
 $d\sigma^2=|\omega|^2$ whose divisor is 
 \[
   \sum_{j=1}^{m} \nu_j q_j+\sum_{l=1}^n p_l.
 \]
 Then the flat front constructed from $(G,\omega)$ is well-defined on 
 \[
   M^2=\overline{M}^2\setminus\{p_1,\dots,p_n,q_1,\dots,q_m\}.
 \]
 Moreover, each $q_j$ is a regular end (see Proposition \ref{prop:regularend}), and 
 since $G$ does not branch at $q_j$, it is a properly embedded end 
(Proposition  3.12 of \cite{KUY2}).  
 Then, since $\ord_{p_l}|\omega|^2=1$ and $p_l$ is a branch point of $G$ with
 multiplicity $1$, \eqref{eq:s_ord}, \eqref{eq:s-order} and 
 \eqref{eq:schwarz} yield that 
 $Q$ has at most a simple pole at $p_l$.  
 So we expect that generically $Q$ has a simple pole at each $p_l$.  
 Hence $\ord_{p_l} |\theta|^2 = -2$
 and Proposition \ref{prop:ordQgem1} gives 
 $m(f,p_l)=1$, implying that $p_l$ is an embedded end.  
 Thus we can expect the existence of
 many higher genus flat fronts with embedded ends.
 If $m=1$, such an $f$ might have genus $k$ with
 $2k+3$ embedded ends (one end from $q_1$ and $2k+2$ ends from the $p_l$).  
 In fact, when $k=1$ such an example is given 
 explicitly in \cite{KUY2}.
 So the next natural question is:
\begin{problem}
 Is there a complete flat front in $H^3$ 
 of genus $k \ge 1$ with $2k+2$ embedded ends?
\end{problem}
\end{remark}
Even when $k=1$, no such examples are known; the problem was raised in 
\cite[Remark~3.17]{KUY2} for $k=1$.  
The following examples are related to this problem.
\begin{example}\label{ex:4k+1}
 We construct flat fronts of genus $k \ge 1$ 
 with $4k+1$ embedded ends, which are canonical generalizations of
 the $5$-ended genus $1$ fronts in \cite{KUY2}.  Choose a polynomial 
 $\varphi(z)$ so that 
 \begin{enumerate}
 \item[(a)]\label{ex:item0} 
           $\varphi(z)=z^{2k}+a_1z^{k-1}+\dots+ a_{k-1}z+a_k$, 
           where $a_1\neq 0$, $a_k\neq 0$, and
 \item[(b)]\label{ex:item1} 
           $\varphi(z)$ has only simple roots, and
 \item[(c)]\label{ex:item2}
           $(z \varphi(z))'=\varphi(z)+z\varphi'(z)$ 
           also has only simple roots.
 \end{enumerate}
 Consider the genus $k$ hyperelliptic Riemann surface defined by 
 $w^2=z\varphi(z)$.  Set 
 \begin{equation}\label{eq:1}
   G=w,\qquad G_*=\frac{h}{w},\quad
   \text{where}\quad
   h=h(z)=\frac{1}{2k+1}{\bigl(2kz\varphi(z)-z^2\varphi'(z)\bigr)}.
 \end{equation}
 Then we have
 \begin{equation}\label{eq:2}
   \frac{dG}{G-G_*}=
   \frac{G\,dG}{G^2-G_*G}
     =\frac{\varphi(z)+z\varphi'(z)}{2(z\varphi(z)-h(z))}=
   \frac{(2k+1)\,dz}{2z}.
 \end{equation} 
 Clearly $G$ is of degree $2k+1$. 
 Since $h/w$ has only simple poles and only at 
 the zeros of $\varphi(z)$, $G_*$ has degree $2k$.  Thus
 $\deg G+\deg G_*=4k+1$.  By \eqref{eq:2}, 
 \begin{equation}\label{eq:3}
   \xi=\exp \int \frac{dG}{G-G_*}=z^{(2k+1)/2} 
 \end{equation}
 has only purely imaginary monodromy.  
 Moreover, \eqref{eq:1} implies 
 \begin{equation}\label{eq:star}
  G-G_*=\frac{1}{w}(w^2-h)=\frac{z}{(2k+1)w}(\varphi+z\varphi'). 
 \end{equation}
 So applying Fact~\ref{fact:g-rep},
 we get a flat front
 \[
   f\colon M^2=\overline M^2
            \setminus\{(z,w)\, ; \, z\bigl(\varphi(z)+z\varphi'(z)\bigr)=0\}
            \longrightarrow H^3.
 \]
 By \eqref{eq:star}, $f$ has exactly $4k+1$ ends, since the
 ends are the zeros of $z\bigl(\varphi(z)+z\varphi'(z)\bigr)$.
 Note that $z=\infty$ is not an end, 
 since $G(\infty)=\infty$ and $G_*(\infty)=0$.  By 
 \eqref{eq:3} and \eqref{eq:data-g}, one canonical $1$-form is 
 $\omega=-z^{-(2k+1)}dG$, and $f$ could have been made using Fact 
 \ref{fact:legendre} as well, 
 with this $\omega$.  As 
 \[
    dG_*=\frac{h'w-hw'}{w^2}dz
     =\frac{h'w^2-hww'}{w^3}dz
     =\frac{h'z\varphi-{h}(z\varphi)'/2}{w^3}dz,
 \]
 we have by \eqref{eq:star}, \eqref{eq:2} and \eqref{eq:data-g} 
 that
 \[
  Q=
  -\frac{2k+1}{2z} \cdot \frac{zh' \varphi -{h}(z\varphi)'/2}{w^2(w^2-h)}dz^2
  =-\frac{(2k+1)^2}{2} \cdot 
  \frac{{h'\varphi}-{(h/z)}(z\varphi)'/2}{z w^2\left(\varphi+z\varphi'\right)}
  dz^2.
 \]
 Here, $dz/w$ has no poles and zeros at the zeros of $w$.
 Since $a_k\ne 0$,  the origin $z=0$ is a pole of order $2$ 
 of $Q$, and thus it is a weakly complete end, by Proposition \ref{prop:order}.  
 On the other hand, $Q$ has simple poles at the other ends.
 However, $\omega$ has zeros at the ends other than $z=0$, by \eqref{eq:2}.  
 Hence the orders of $|\theta|^2$ are less than $-1$ at such ends,
 and then they are weakly complete.
 Any branch point of $G$ in $\overline{M}^2$ 
 must be a zero of $(z \varphi)'$, and thus must be an end of $f$, 
 so $G$ has no branch points on $M^2 \setminus \{ z=\infty \}$.
 Moreover $a_1a_k\ne 0$ implies that $G_*$ does not 
 branch at $z=\infty$ and $z=0$.  
 Thus $f$ has no branch points anywhere on $M^2$, by Theorem 2.9 of \cite{KUY2}.  
 Since $f$ is weakly complete and of finite type, 
 there exist infinitely many complete fronts in the 
 parallel family of $f$, by Theorem \ref{thm:fff->f_t}.  
 By equality of the Osserman-type inequality (\cite[Theorem~3.13]{KUY2}), 
 all the ends of $f$ are properly embedded.  
\begin{figure}
 \begin{center}
  \includegraphics[width=5cm]{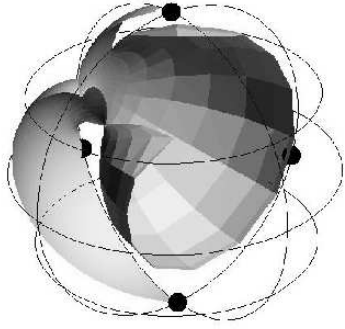}
 \end{center}
 \caption{
 A genus $2$ flat front with $10$ embedded ends, coming from a 
 variant of the approach in Example~\ref{ex:4k+1}.  
 It is defined on the 
 closed Riemann surface 
 $\overline{M}^2 = 
 \{ (z,w) \in (\C \cup \{ \infty \})^2 \, ; 
 \, w^2 = z (z^2-1)(z^2-9/4) \}$ 
 with $10$ points removed, and with 
 data $G=w$ and $G_* = (5w-2 z (dw/dz))/5$.  
 Because $G$ and $G_*$ have no common branch points, the surface is a
 front with embedded ends.  
 The portion shown here is the image of one sheet in $\overline{M}^2$ over 
 the quadrant 
 $\{ z \in \C \cup \{\infty\} \,; \, \Re z \geq 0, \Im z \geq 0 \}$ 
 in the $z$-plane, and is one-eighth of the full 
 surface.  
 The planar boundary curves of this portion lie in planes of reflective 
 symmetry of the full surface.  
 The four ends appearing in the boundary of 
 this portion are marked with black dots.  
}\label{fig:genus2}
\end{figure}

 Finally, we shall show that
 $\varphi(z)=z^{2k}-2c z^{k-1}-1$
 satisfies the above conditions (b), (c) 
 for a suitable $c\in \R$.  When $k=1$ and $c=0$, $f$ is just the example in 
 \cite[Example~4.6]{KUY2}.
 When $k\ge 2$, we set $c=k/(k-1)$. 
 Then 
 \[
   \varphi'(z)=2 k z^{k-2} \bigl( z^{k+1}-1 \bigr).
 \]
 Let $p \neq 0$ be a zero of $\varphi'(z)$, 
 so $p^{k+1}=1$  and $p^2\varphi(p)=1-p^2-2k/(k-1)$.
 If $\varphi(p)=0$, then $p^2=(1+k)/(1-k)$ 
 and $|p| \neq 1$.  So $\varphi(z)$ has no double roots.  
 Now 
 \[
    (z\varphi)''=2kz^{k-2}((2k+1)z^{k+1}-k).
 \]
 Let $q \ne 0$ be a zero of $(z\varphi)''$.  
 Then we have $q^{k+1}=k/(2k+1)$.
 Hence $|q|^2=\bigl(k/(2k+1)\bigr)^{2/(k+1)}$ is not a
 rational number.  On the other hand, 
 \[
   q^2 (z\varphi)'|_{z=q}=(2k+1)\left(\frac{k}{2k+1} \right)^2-\frac{2 k^2}{k-1} 
   \cdot \frac{k}{2k+1}
         -q^2.  
 \]
 If $(z\varphi)'|_{z=q}=0$, we have 
 $q^2=k^2(1+k)/((2k+1)(1-k))\in \Q$, a contradiction. 
 Thus both $\varphi$ and $(z\varphi)'$ have only simple roots, 
 and $\varphi(z)$ satisfies conditions (b), (c).  
\end{example}

\section{p-fronts}
\label{sec:p-front}
Now we consider flat p-fronts, starting with a proof of Theorem \ref{thm:B}.  
\begingroup
\setcounter{introtheorem}{1}
\begin{introtheorem}
Let $f\colon M^2 \to H^3$ be a flat p-front, then $M^2$ 
is orientable. 
\end{introtheorem}
\endgroup
As noted in the introduction, the other space forms 
$\R^3$ and $S^3$ admit flat M\"obius bands (see G\'alvez and Mira \cite{GM}).  
So Theorem \ref{thm:B} is special to $H^3$.  Kitagawa \cite{Kit} proved 
the orientability of compact flat surfaces in $S^3$. 
\begin{proof}[Proof of Theorem~\ref{thm:B}]
 As $f$ is a p-front, for any $p \in M^2$ 
 there exists a neighborhood $U_p \subset M^2 $ of $p$ such that 
 $f|_{U_p}$ is a front.  We may assume that $U_p$ is simply connected. 
 Then as noted in Section~\ref{sec:prelim},
 there exists a unique complex structure on $U_p$
 such that both hyperbolic Gauss maps $G$ and $G_*$ are meromorphic.  

 Since $f|_{U_p}$ is a front, at least one of $G$ or $G_*$ is
 not branched at $p$, by Theorem 2.9 of \cite{KUY2}.  
 Without loss of generality, 
 we may assume $G$ and $G_*$ are finite at $p$.  
 Choosing $U_p$ sufficiently small, we have a local complex 
 coordinate $z=G|_{U_p}$ or $z=G_*|_{U_p}$ on $U_p$ at each point $p\in M^2$. 
 Since $G\circ G^{-1}$ and $G_*\circ G^{-1}_*$
 are identity maps and either $G\circ G^{-1}_*$ or
 $G_*\circ G^{-1}$ is well-defined and holomorphic 
 on $U_p$ for each $p\in M^2$,
 the transition function on $U_p\cap U_q$
 for two distinct points $p$ and $q$ is always 
 holomorphic.
 So we can extend this local complex structure on $U_p$ to all of $M^2$.
\end{proof}
When we consider a flat p-front,
we always regard $M^2$ as a Riemann surface with the complex structure 
given in the proof of Theorem~\ref{thm:B}.
Note that co-orientability is defined in the introduction.
\begin{theorem}\label{thm:double-lift}
 Let $f\colon M^2 \to H^3$ be a non-co-orientable flat p-front 
 with universal cover $\pi\colon{}\widetilde M^2 \to M^2$. 
 Then there exists a Legendrian immersion
 \[
    \E_f:\widetilde{M}^2\longrightarrow \SL(2,\C)
 \]
 and a covering transformation 
 $\tau\colon{}\widetilde{M}^2\to \widetilde{M}^2$
 such that 
 \[
       \E_f\circ \tau= 
             \E_f \begin{pmatrix}
                     0 & \iii \\ \iii & 0 
                  \end{pmatrix}\,
                  (= \E^{\natural}_f). 
 \]
\end{theorem}
We call this $\E_f$ the {\it adjusted lift\/} of the p-front $f$. 
Conversely, if a lift $\E_f$ satisfies 
$\E_f\circ \tau=  \E^{\natural}_f$ for some covering transformation $\tau$, 
$f=\E^{}_f\E_f^*=\E^{\natural}_f(\E^{\natural}_f)^*$
is non-co-orientable. 
\begin{proof}[Proof of Theorem~\ref{thm:double-lift}]
 Take a holomorphic Legendrian lift 
 $\E_0\colon{}\widetilde{M}^2\to \SL(2,\C)$ of 
 $f$; it is determined up to right-multiplication by a matrix in $\SU(2)$.  
 We now change $\E_0$ to an adjusted lift.  
 The flat front $\tilde f=\E^{}_0\E_0^*\colon{}\widetilde{M}^2\to H^3$ 
 satisfies $f\circ \pi=\tilde f$. 
 The unit normal vector of 
 $\tilde f$ is $\tilde \nu=\E^{}_0e_3\E_0^*$, 
 by \eqref{eq:front-normal}.  
 If $p_1,p_2\in \pi^{-1}(p)$,
 then $\tilde\nu(p_1) = \pm \tilde\nu(p_2)$.  
 So there exists a (unique) representation
 $\delta \colon{}\pi_1(M^2)\to \{\pm 1\}$ such that
 \[
   \tilde \nu \circ T=\delta (T) \tilde \nu \qquad (T\in \pi_1(M^2)),
 \]
 where the fundamental group $\pi_1(M^2)$ is identified with the 
 covering transformation group.  
 Since $f$ is non-co-orientable, $\delta$ is surjective.
 Letting $\check M^2 = \widetilde{M}^2/\Ker\delta$ with associated 
 nontrivial covering involution $\sigma$ on $\check M^2$, 
 $\check \nu=\E^{}_0e_3\E_0^*$ is single-valued on $\check M^2$, 
 so $\check f = \E^{}_0\E_0^* \colon{}\check M^2\to H^3$ is a flat front
 because its unit normal vector $\check \nu$ is well-defined, and 
 $\check \nu\circ \sigma=-\check \nu$.
 Since $\widetilde{M}^2$ is universal, there exists a lift
 $\tau\colon\widetilde{M}^2\to \widetilde{M}^2$ of $\sigma$.
 As $\tau$ is holomorphic and 
 $\E^{}_0\E^*_0=(\E^{}_0\circ\tau)(\E^{}_0\circ\tau)^*$, 
 there exists an $\SU(2)$-matrix
 \[
   R=\begin{pmatrix} p & -\bar q \\ q &\hphantom{-}\bar p \end{pmatrix}
       \qquad (|p|^2+|q|^2=1)
 \]
 such that 
 $\E^{}_0\circ \tau=\E_0 R$.
 Since $\check \nu\circ \sigma=-\check \nu$, we have 
 $\E^{}_0e_3\E_0^*= - \E_0Re_3R^*\E_0^*$, so 
 $Re_3R^*= - e_3$, implying $p=0$.  Thus 
\[
  \E^{}_f=\E^{}_0 a,\qquad
 \text{where}\quad
 a=\begin{pmatrix} \sqrt{i/q} & 0 \\ 0 & \overline{\sqrt{i/q}} \end{pmatrix},
\]
 satisfies $f=\E^{}_f\E^*_f$ and
 $\E^{}_f\circ \tau=\E^{\natural}_f$.
\end{proof}
\begin{corollary}\label{cor:double-lift}
 Let $f\colon M^2 \to H^3$ be a non-co-orientable flat p-front. 
 Then there is a double cover $\check \pi\colon \check M^2 \to M^2$ 
 such that $\check f = f \circ \check \pi\colon \check M^2 \to H^3$ 
 is a front.
 Moreover, there exists a covering transformation 
 $\tau \colon \check M^2 \to \check M^2$ with $\check f \circ \tau = \check f$ 
 satisfying 
 \[
  \check G \circ \tau = \check G_*, \quad
  \check G_* \circ \tau = \check G, \quad 
  \check \omega \circ \tau =\check \theta, \quad 
  \check \theta \circ \tau = \check \omega \quad\text{and}\quad
  \check Q \circ \tau = \check Q, 
 \]
 where $\check G = G \circ {\check \pi}$ and $\check G_*$, $\check \omega$, 
 $\check \theta$, $\check Q$ are defined similarly.  
\end{corollary}
Thus, for a p-front, the Hopf differential $Q$ and 
$ds^2_{1,1}=|\omega|^2+|\theta|^2$ are well-defined.
\begin{definition}\label{def:complete-p-front}
 A non-co-orientable flat p-front $f$ is called 
 {\it complete\/} (resp. {\it weakly complete, finite type}) if 
 its double cover 
 $\check f$ as in Corollary \ref{cor:double-lift} is 
 complete (resp. weakly complete, finite type).  
\end{definition}

\begin{figure}
\begin{center}
   \begin{tabular}{c@{\hspace{1cm}}c}
        \includegraphics[width=3cm]{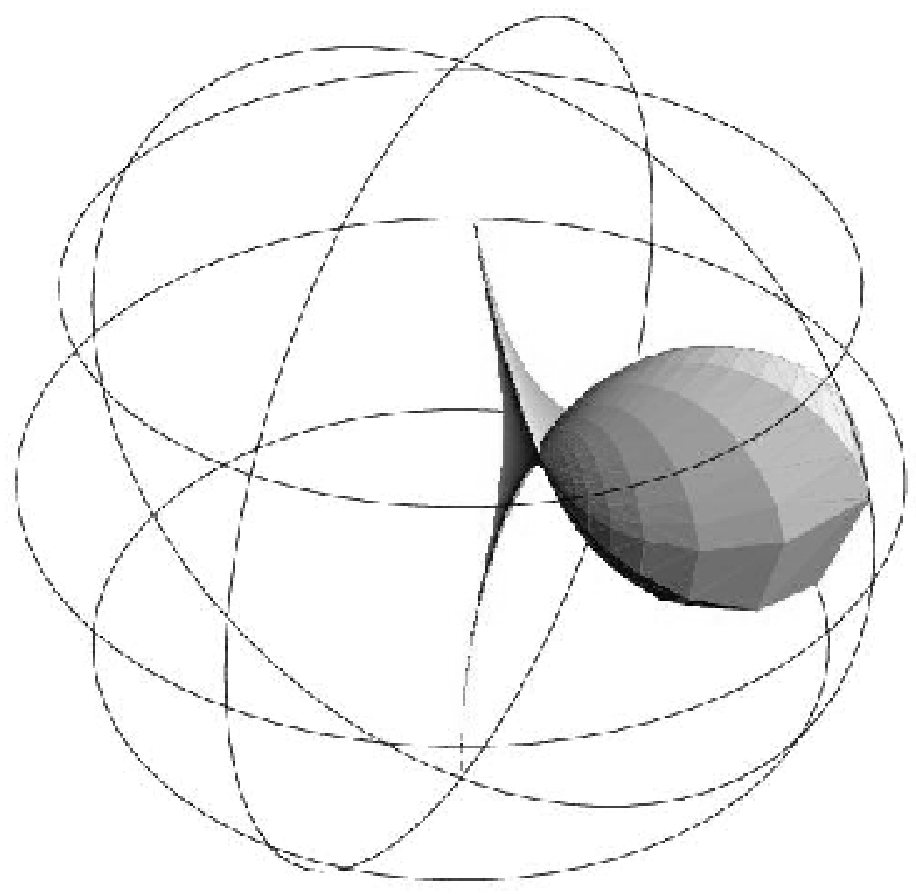}&
        \includegraphics[width=3cm]{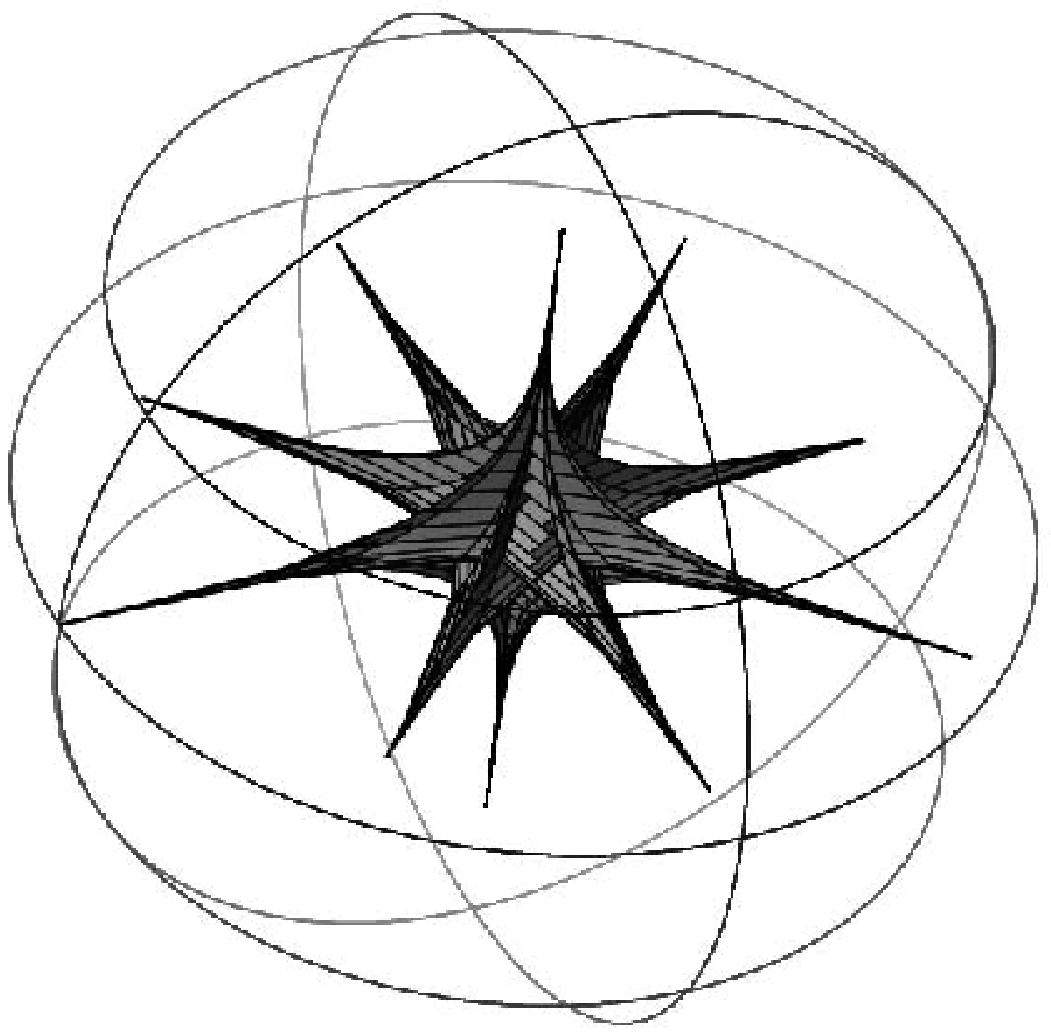}
  \end{tabular}
\end{center}
\caption{%
 The p-front in Example~\ref{ex:peach} that is not globally 
 a caustic on the left, and the caustic with dihedral cross $\Z^2$
 symmetry for the fronts produced by $G=z^3$ and $G_*=z^{-5}$ and 
 $M^2=\C \setminus \{ z \,;\, z^8=1 \}$ on the right.
}
\label{fig:not-caustic}
\end{figure}

\begin{proposition}\label{inserted-at-last-minute-proposition}
Let $f\colon{}M^2\to H^3$ be a flat p-front.  The following hold{\rm :}
\begin{enumerate}
 \item\label{last:1} 
       If $f$ is complete, then it is 
       weakly complete and of finite type.
 \item\label{last:2}
       If $f$ is weakly complete and of finite type, there exists a 
       compact Riemann surface $\overline{M}^2$ such that
       $M^2$ is biholomorphic to 
       $\overline{M}^2\setminus\{p_1,\dots,p_n\}$,
       for finitely many distinct points $p_1,\dots,p_n$ in 
       $\overline{M}^2$.
       Moreover, $ds^2_{1,1}$ has finite order at $p_1,\dots,p_n$ 
       and $Q$ extends to a meromorphic 
       $2$-differential on $\overline{M}^2$.
\end{enumerate}
\end{proposition}
\begin{proof}
 \ref{last:1} follows immediately from Proposition~\ref{prop:f->finite}. 
 To show \ref{last:2}, 
 $M^2$ is biholomorphic to $\overline{M}^2\setminus\{p_1,\dots,p_n\}$ by
 Huber's theorem, 
 just as in the proof of Proposition~\ref{prop:finite->}.  
 We will call the points $\{p_1,\dots,p_n\}$ the {\em ends\/} of
 $f$ (Definition~\ref{def:puncture}).
 Suppose  an 
 end $p_j$ is co-orientable, 
 that is, the restriction of $f$ to a sufficiently small punctured
 neighborhood of $p_j$ is co-orientable 
 (see Definition~\ref{def:co-orientable-end}).
 Then Proposition~\ref{prop:finite->} implies
 $\omega$, $\theta$ have finite order at $p_j$, and $Q$ is meromorphic
 there.  
 If $p_j$ is not co-orientable, 
 we take a punctured neighborhood $U_j^*$ of $p_j$, 
 and let $\check U_j^*$ be its double cover.  
 Lifting $f|_{U_j^*}$ to $\check f : \check U_j^* \to H^3$, $\check f$
 is now a weakly complete co-orientable end of finite type.  
 Its canonical $1$-forms $\check \omega$, $\check \theta$ are of finite
 order, and 
 $d\check s^2_{1,1}=|\check \omega|^2+|\check \theta|^2$ is complete and
 of finite order, and its Hopf differential $\check Q$ extends
 meromorphically across the puncture of $\check U_j^*$, by
 Proposition~\ref{prop:finite->}.
 Noting that $d\check s^2_{1,1}$ is actually well defined on $U_j^*$ and 
 equals $ds^2_{1,1}$ there, and that $\check Q$ projects to the Hopf
 differential of $f|_{U_j^*}$ on $U_j^*$, the proof is completed.  
\end{proof}
To observe the behavior of ends of $p$-fronts, 
we review properties of regular ends of fronts:
\begin{definition}\label{def:puncture}
A p-front $f\colon{}M^2\to H^3$ is said to be of {\it finite topology\/} if 
there exists a compact Riemann surface $\overline{M}^2$ such that
$M^2$ is homeomorphic to
$\overline{M}^2\setminus\{p_1,\dots,p_n\}$.
Such points $p_1,\dots,p_n$ are called the {\it ends\/} of $f$.

 An end of a weakly complete flat p-front of finite topology is said to be 
 of {\it puncture-type\/} if it is biholomorphic to a punctured disk, 
and to be {\it annular\/} if it is biholomorphic to an annulus 
$\{ z \, ; \, r_1 < |z| < r_2\}$, $0<r_1<r_2<\infty$.  

 A puncture-type end $p$ of a weakly complete flat front is said to be
 \emph{regular} if both $G$, $G_*$ have at most poles at $p$,  
 and to be \emph{irregular} otherwise. 

 For a weakly complete flat p-front, a puncture-type end $p$ 
 is said to be {\it regular\/} if the corresponding end 
 $\check p$ of $\check f$ is regular, 
 and to be \emph{irregular} otherwise.  
\end{definition}
\begin{proposition}\label{prop:regularend}
 For an end $p$ of a flat front of finite type, 
 the following conditions are equivalent{\rm : }
 \begin{enumerate}
  \item The end $p$ is regular. 
  \item One of $G$, $G_*$ has at most a pole at $p$. 
  \item The Hopf differential $Q$ has at most a pole of  
        order $2$ at $p$. 
 \end{enumerate}  
\end{proposition}
\begin{proof}
This was proven in \cite[Theorem 4]{GMM} for complete ends.  
Reading that proof carefully, one sees that it applies to 
flat fronts of finite type as well.  
\end{proof}
\begin{proposition}
 Let 
 $f \colon M^2 = \overline M^2 \setminus \{p_1, \dots , p_n \} \to H^3$ 
 be a weakly complete flat front 
 whose ends are all regular. 
 Then $f$ is of finite type if and only if the Hopf differential has at
 most a pole of order $2$ at each end.  
\end{proposition}
\begin{proof}
 One direction follows from Proposition \ref{prop:regularend}, 
 so we prove the other direction here.  Choose an end $p_j$. 
 By \eqref{eq:s-order}, $\omega$ and $\theta$ have 
 finite order at $p_j$ if and only if $s(\omega)$ and $s(\theta)$ have at 
 most poles of order $2$.  
 Since $G$, $G_*$ are meromorphic at $p_j$,
 $S(G)$ and $S(G_*)$ have at most poles of order $2$.
 Then \eqref{eq:schwarz} implies that $\ord_{p_j} Q \geq -2$ 
 if and only if $s(\omega)$ and $s(\theta)$ have at most poles of order $2$.  
\end{proof}

\begin{definition}\label{def:co-orientable-end}
 For a weakly complete flat p-front
 $f\colon{}\overline{M}^2\setminus \{p_1,\dots,p_n\}\to H^3$, 
 an end $p_j$ is called {\it co-orientable\/} if
 the restriction of $f$ to a sufficiently small punctured
 neighborhood of $p_j$ is co-orientable, and is otherwise 
 called {\it non-co-orientable}.
\end{definition}

A co-orientable regular end $p$ of a flat p-front $f$
can be considered as a regular end of a flat front, and 
we have already defined its multiplicity $m(f,p)$ 
(Definition~\ref{def:multiplicity}).  
We also now define the multiplicity of a non-co-orientable end.
Let 
\[
f\colon{}D^*(\varepsilon)=\{z\,;\,0<|z|<\varepsilon\} \longrightarrow H^3
\]
be a non-co-orientable end at $p$.
Then there exists a lift (as a p-front)
$\check f\colon{}\check D^*(\varepsilon)\to H^3$
of $f$, where $\check D^*(\varepsilon)$ is the double covering of 
$D^*(\varepsilon)$.
We set
\begin{equation}\label{eq:multiplicity-non-co-orientable}
   m(f,p) = m(\check f , \check p)/2
\end{equation}
and call it  the \emph{multiplicity} of the end $p$ of $f$. 
Thus we have defined the multiplicity of 
any regular end of a weakly complete p-front, taking its value in 
$\frac{1}{2}\Z$. 
Let $f\colon{}M^2\to H^3$ be a weakly complete flat front with regular 
ends.  
If $f$ is co-orientable, the hyperbolic Gauss maps $G$, $G_*$ of $f$
are single-valued on $M^2$ and we set
\[
    \deg \G_{f} = \deg G + \deg G_{*},
\]
which we call the {\it total degree} of the Gauss maps of $f$.
When $f$ is non-co-orientable, there exists a
lift $\check f\colon{}\check M^2\to H^3$ of $f$, 
where $\check M^2$ double covers $M^2$, and we set 
\[
     \deg \G_{f} = (\deg \mathcal G_{\check f})/2 \; .  
\]
It follows from \cite[Theorem 3.13]{KUY2} that 
\begin{equation}\label{p-Osserman}
   \deg {\mathcal G}_f \ge \# \{\text{\rm non-co-orientable ends}\}/2 + 
   \# \{\text{\rm co-orientable ends}\}
\end{equation}
 for a weakly complete flat p-front of finite type 
whose ends are all regular. 
When $f$ is complete, equality implies all ends are properly embedded.
We note that 
{\it complete ends are automatically co-orientable}: 
Suppose a complete end $p$ of a p-front $f\colon{}M^2\to H^3$
is non-co-orientable.
Take the double cover $\check M^2$, the lift 
$\check f\colon{}\check M^2\to H^3$  
and a complex coordinate $z$ of $\check M^2$ around the (unique) lift $\check p$
of $p \in \overline{M}^2$ with $z(\check p)=0$.  
Since $\check p$ is also a complete end, $\check f$ is of finite type
and the canonical forms are written as
\[
    \check \omega = cz^{\mu}\bigl(1+o(1)\bigr)\,dz
    \qquad\text{and}\qquad
    \check \theta = c_* z^{\mu_*}\bigl(1+o(1)\bigr)\,dz,
\]
where $c,c_*$ are non-zero constants and $\mu,\mu_*\in\R$.
Here, by Corollary~\ref{cor:double-lift}, 
we have $|c|=|c_*|$ and $\mu=\mu_*$.
Hence we have 
$|\check \rho| =|\check \theta/\check \omega|=1+o(1)$, 
which implies that 
the singular set of $\check f$ accumulates at $\check p$, contradicting 
that $\check p$ is a complete end.  

\begin{example}[Weakly complete p-fronts with $3$ ends]
 \label{ex:3-noid} This example is of interest, because it is neither 
 a front nor globally a caustic (see Remark \ref{noncaustic-nonfront}).  Set
\[
  G = \frac{z^2+\tfrac{z}{b}}{z+b} 
 \quad \text{and}\quad
  G_* = \frac{z^2-\tfrac{z}{b}}{b-z}
  \qquad (b\in\R\setminus\{0,\pm 1\})
\]
 which are defined on the Riemann surface 
 $\check M^2=\C \setminus \{0, \pm 1 \}$.  
 Then Fact~\ref{fact:g-rep} 
 gives a flat front $\check f_b\colon{}\check{M}^2\to H^3$
 with hyperbolic Gauss maps $G$ and $G_*$.  
 We also have 
\begin{align*}
 \xi    &=c \cdot \exp \int \frac{dG}{G-G_*} 
    = c \frac{\sqrt{z}}{z+b} (z-1)^{(1-b)/2} (z+1)^{(1+b)/2}, \\
 \omega &=-c^{-2}(z^2+2bz+1)z^{-1} 
             (z-1)^{b-1}(z+1)^{-b-1}\,dz, \\
 \theta &= -c^2(z^2-2bz+1)z^{-1}
                 (z-1)^{-b-1}(z+1)^{b-1}\,dz/4. 
\end{align*} 
 Now we set $c=\sqrt{2}$.
 Then $\omega\circ\tau=\theta$ holds,
 where $\tau\colon{}\check M^2\ni z\mapsto -z\in\check{M}^2$.  
 Then $\E^{}_{f_b}\circ\tau=\E^{\natural}_{f_b}$, and hence 
 we have a well-defined flat p-front 
 \[
    f_b\colon{}M^2 = (\check M^2/\sim) \longrightarrow H^3
 \]
 where $z_1\sim z_2$ if and only if $z_2=\pm z_1$.
 The three ends of $f_b$ are at $z=0$ and $z=\infty$ and $z = \pm 1$. 
 The ends at $z=0$, $z=\infty$ are non-co-orientable, and the end at
 $z=\pm 1$ is co-orientable. 
 $f_b$ satisfies equality in \eqref{p-Osserman}.
\end{example}

\section{Caustics}
\label{sec:caustic}
Roitman~\cite{R} showed that the caustic $C_f$ 
(or focal surface) of a flat surface $f$ is also flat, and  
gave a representation for $C_f$ in terms of $f$.  
In \cite[Section~5]{KRSUY}, Roitman's representation is described 
in the terminology below: 
Let $f\colon M^2\to H^3$ be a flat front with hyperbolic 
Gauss maps $G,G_*$. 
Let $q_1,\dots,q_m\in M^2$ be the umbilic points of $f$.
Then we can choose a single-valued square root 
\[
           \beta=\sqrt{dG/dG_*}
\]
defined on the universal cover $\McUniv$ of 
$\Mc=M^2 \setminus \{q_1,\dots,q_m\}$.
The caustic $C_f$ is 
\begin{equation*}
 C_f=\Ec {\Ec}^* \colon \Mc 
     \longrightarrow H^3 \; , 
\end{equation*}
\begin{align}\label{eq:Ec}
 \Ec 
  = \frac{\sqrt{i}}{\sqrt{2\beta(G-G_*)}}
 \begin{pmatrix}
  G + \beta G_* & \iii (G - \beta G_*) \\
  1 + \beta & \iii (1 - \beta) 
 \end{pmatrix}
 \begin{pmatrix}
  \sqrt{i} & 0 \\
  0 & 1/\sqrt{i} 
 \end{pmatrix} \in \PSL(2,\C) \; , 
\end{align}
where $\sqrt{i}=e^{\pi i/4}$.  Note that $\Ec \in \PSL(2,\C)$ because of the 
sign ambiguity of $\sqrt{2 \beta (G-G_*)}$.  
The $\SU(2)$-matrix $\diag(\sqrt{i},1/\sqrt{i})$
in Equation \eqref{eq:Ec} is not essential, 
and is included only so that $\Ec$ changes to $\Ec^\natural$ when
$\beta$ changes to $-\beta$, i.e. so that $\Ec$ becomes an adjusted lift.  
The hyperbolic Gauss maps $\Gc, \Gc_*$ and the canonical forms $\omegac$, 
$\thetac$ of the caustic $C_f$ are 
\begin{equation}\label{hGm_c}
 (\Gc, \Gc_*)  = 
  \left( \frac{G + \beta G_*}{1+\beta},
       \frac{G - \beta G_*}{1-\beta} \right), 
\end{equation}
\begin{align}
 \omegac 
   & = \frac{1}{4}\left\{
       \frac{2 ( \beta+1 )^2}{G-G_*}dG_*
            -d \log \left(\frac{dG}{dG_*}\right)
           \right\}, \label{omega_c}\\  
 \thetac & = 
        \frac{1}{4}\left\{
        \frac{2 ( \beta-1 )^2}{G-G_*}dG_*
            -d \log \left(\frac{dG}{dG_*}\right)
        \right\}. \label{theta_c} 
\end{align}
These two forms can be also expressed using 
$Q (= \omega \theta)$, $\rho (=\theta / \omega)$
of the original front $f$, as can the Hopf differential $\Qc 
= \omegac \thetac$ of $C_f$: 
\begin{equation}\label{omega_c-theta_c}
 \omegac = i \sqrt{Q} + \frac{1}{4} d \log \rho, \quad
 \thetac = - i \sqrt{Q} + \frac{1}{4} d \log \rho, \quad
 \Qc = 
 Q + \left(\frac{d \log \rho}{4} \right)^2 , 
\end{equation}  
where the sign of $\sqrt{Q}$ is chosen so that \eqref{omega_c-theta_c} is 
compatible with \eqref{omega_c} and \eqref{theta_c}.  
\begin{proposition}\label{prop:caustic-front}
 The caustic of a flat front is a p-front.  
\end{proposition}
\begin{proof}
 By Theorem 2.9 of \cite{KUY2}, it 
 suffices to prove that $\omegac$ and $\thetac$ have no common zeros. 
 Let $p$ be an arbitrary point on the caustic. 
 Since $p$ is not an umbilic point on the original front, i.e., 
 $Q(p) \ne 0$, it follows from \eqref{omega_c-theta_c} that  
 at least one of $\omegac(p), \thetac(p)$ is not zero.
\end{proof}
\begin{example}[Caustics of rotational flat fronts]\label{example6pt2}
 Since the horosphere is totally umbilic, it has no caustic.  For other 
 rotational examples: 
\begin{enumerate}
 \item the caustic of a hyperbolic cylinder is a hyperbolic line, 
 \item the caustic of an hourglass is also a hyperbolic line,
 \item the caustic of a snowman is a hyperbolic cylinder, 
\end{enumerate} 
 where hourglasses and snowmen are fronts of revolution with single conical 
 singularities and cuspidal edge singularities, respectively 
 (see \cite[Example 4.1]{KUY2}).  In fact, snowmen are 
 characterized by their caustics, as Proposition \ref{prop:cylinders} shows:  
\end{example}

\begin{proposition}\label{prop:cylinders}
 Let $f$ be a flat front with caustic $\Cf$.  
 Then $f$ is locally a snowman if and only if $\Cf$ is an open
 submanifold of a hyperbolic cylinder.  
\end{proposition}
\begin{proof}
 As seen in Example \ref{example6pt2}, a flat front of revolution has a 
 hyperbolic cylinder as its caustic if and only if it has a 
 nonempty cuspidal edge set, so we need only show that if $\Cf$ 
 is a hyperbolic cylinder, then $f$ is a surface of revolution.  

 Assume $\Cf\colon{}M^2\to H^3$ is a hyperbolic cylinder, so we may
 assume 
 $M^2$ is an open subset of $\C \setminus \{ 0 \}$
 and 
 \[
     \thetac = \frac{c^2}{4z}\,dz ,\qquad
     \omegac = \frac{1}{c^2z}\,dz,
 \]
 where $c\in \R\setminus\{0\}$ (see \cite[Example 4.1]{KUY2}).
 Then by \eqref{omega_c-theta_c}, the Hopf differential $Q$ of $f$
 is
 \[
    Q = -\frac{1}{4}(\thetac-\omegac)^2 =
        -\frac{a^2}{4}\frac{dz^2}{z^2}
     \qquad \left(a=\frac{c^2}{4}-\frac{1}{c^2}\right),
 \]
 and the function $\rho$ satisfies
 \[
    d\log \rho = \frac{2b}{z}\,dz,
    \qquad\text{and then}\qquad
     \rho = -k z^{2b}
    \qquad
    \left(b=\frac{1}{c^2}+\frac{c^2}{4}\right),
 \]
 where $k\neq 0$ is a constant.
 Thus, we compute $\theta=\sqrt{Q\rho}$ and $\omega=Q/\theta$ as
 \[
    \theta =  \frac{\sqrt{k}}{2}az^{b-1}\,dz, \qquad
    \omega =  \frac{1}{2\sqrt{k}}az^{-b-1}\,dz.
 \]
 Hence $f$ is a part of a surface of revolution.
\end{proof}

One can check that the caustic of a peach front (Example~\ref{ex:peach}) is a
horosphere \cite{R}.  Peach fronts are also characterized by their caustics:
\begin{proposition}\label{prop:peach}
 A flat front is locally a peach front 
 if and only if its caustic is totally umbilic, i.e., is 
 locally an open submanifold of a horosphere.  In particular, a 
 peach front is the only flat front whose caustic is the horosphere. 
\end{proposition}
\begin{proof}
 Assume the caustic is totally umbilic, that is, $\Qc =0$. 
 Taking the dual lift if necessary, we may assume $\Gc$ is constant, and then 
 applying a rigid motion of $H^3$ if necessary, we may assume $\Gc=0$. 
 It follows from \eqref{hGm_c} that 
 $G^2/G_*^2 = dG / dG_*$, so $1/G_* - 1/G$ is a constant. 
 By a suitable motion of $H^3$, we can change $G$ and $G_*$ to 
 $1/G_*$ and $1/G$, and then $G_* = G + \text{constant}$.  If $G$ branches, 
 then $G_*$ also branches, so the flat surface $f$ produced by $G$ and $G_*$ 
 would not be a front, by Theorem 2.9 of 
 \cite{KUY2}.  Hence we may assume $G=z$ is a coordinate for the front $f$, 
 proving the assertion.
\end{proof}

\begin{remark}
 Propositions \ref{prop:cylinders} and 
 \ref{prop:peach} show that if the caustic $\Cf$ is complete without
 singularities 
 (so is a cylinder or horosphere, see \cite[Theorem 3]{GMM}), then the
 original front $f$ must be a snowman or peach front.  
\end{remark}

\begin{theorem}\label{thm:pf-is-l-c}
 Any flat p-front is locally the caustic 
 of some flat front.  
\end{theorem}  
\begin{proof}
 For $f \colon M^2 \to H^3$ a flat p-front, take any $p \in M^2$. 
 Take a simply connected neighborhood $U$ of $p$, so the canonical forms 
 $\omegac$, $\thetac$ of $f$ are single-valued on $U$.  Set 
\begin{equation*}
 Q_{s} = - \tfrac{1}{4}(e^{\iii s} \omegac 
  - e^{- \iii s} \thetac)^2 \qquad (s\in \R).
\end{equation*}
 Then $Q_{s}$ is a 
 well-defined holomorphic $2$-differential on $U$. 
 If $Q_{s}(p) =0$ for all $s$, 
 then $\omega(p)=\theta(p)=0$, a contradiction, so 
 we can choose a real number $s_0$ so that $Q_{s_0}(p) \ne 0$. 
 Determine  $\omega$ and $\theta$ by
\[
 \omega  \theta = Q_{s_0}, \quad
 \frac{\theta}{\omega}=\rho,\quad
 \text{where}\quad
 \rho=\rho(z)=\exp \left(\int_{z_0}^z 2 
  (e^{\iii s_0}\omegac + e^{-\iii s_0}\thetac)\right) (\ne 0). 
\]
 These $\omega$, $\theta$ yield a flat front $F=\E \E^*$ 
 by solving \eqref{eq:ode}, 
 and the caustic of $F$ is $f$, up to a rigid motion of $H^3$.  
\end{proof} 

\begin{remark}
 The above proof implies that for a given p-front $f$,
 there is locally a two parameter family of fronts $F_{t,s}$ whose caustics 
 are $f$.
 One of these parameters is the signed-distance $t$ of parallel fronts
 and the other is the $s$ in the proof above.
\end{remark}

\begin{remark}\label{noncaustic-nonfront}
 The ``locally'' condition is necessary in Theorem~\ref{thm:pf-is-l-c}, as 
 there are flat p-fronts that are not caustics globally. 
 In fact, the $f_b\colon  M^2 \to H^3$ defined in
 Example \ref{ex:3-noid} is not a caustic 
 over $M^2$ if $b \not\in \frac{1}{2}\Z$: 
 Suppose, by way of contradiction, that $f_b$ is a caustic of 
 some $f_{\orig} \colon M' \to H^3$. 
 By \eqref{omega_c-theta_c}, we have  
 \[
   4 \sqrt{Q_{\orig}} = i \Bigl(\frac{z-1}{z+1}\Bigr)^b 
   \left\{\frac{z^2+2bz+1}{z(z^2-1)} - 
   \frac{z^2-2bz+1}{z(z^2-1)}\Bigl(\frac{z+1}{z-1}\Bigr)^{2b}  \right\}dz.
 \]
 It follows that $Q_{\orig}$ is not well-defined at $z=0$ if 
 $b \not\in \frac{1}{2}\Z$, contradicting that $f_{\orig}$ is 
 well-defined.  Hence the class of p-fronts is strictly 
 larger than the class of flat fronts and their caustics.  
\end{remark}

\section{Ends of Caustics}
\label{sec:end}
First, we recall the properties of regular ends from \cite{GMM} and \cite{KRUY}.  
Recall that if the meromorphic $2$-differential $Q$ on $\overline{M}^2$ expands as 
\[
   Q = z^k\left\{q_0+q_1 z + o(z)\right\}\,dz^2
   \qquad (q_0\neq 0)
\]
in a complex coordinate $z$ at a point $p \in \overline{M}^2$ with $z(p)=0$, the 
integer $k$ is called the {\em order\/} of $Q$ at $p$ and is denoted $\ord_p Q$.
If $k=-2$, the number $q_0$ is independent of the choice 
of coordinate system.  
We call $q_0$ the {\em top term coefficient\/} of $Q$ at $p$.

\begin{definition}\label{def:end-behavior}
 A weakly complete regular end $p$ of finite type is 
 \emph{cylindrical} if $\omega$ and $\theta$ have the same order 
 at $p$, i.e., $\ord_p |\omega|^2 = \ord_p |\theta|^2$.  
 A complete regular end is 
 \begin{enumerate}
  \item {\it horospherical\/} if $\ord_p Q \ge -1$,
  \item of {\it hourglass-type\/} if it is non-cylindrical with 
        $\ord_p Q = -2$ and positive top term coefficient of $Q$, 
  \item of {\it snowman-type\/} if it is non-cylindrical with $\ord_p Q = -2$ and 
        negative top term coefficient of $Q$. 
 \end{enumerate}
\end{definition}

The hyperbolic cylinder, horosphere, hourglass and snowman
(\cite[Example~4.1]{KUY2}) 
are typical examples of such ends, respectively.

\begin{fact}[\cite{GMM,KRUY}]
 A complete end of finite type is an asymptotic covering of a hyperbolic 
 cylinder, horosphere, hourglass or snowman when the end is
 cylindrical, horospherical, of hourglass-type or of snowman-type, respectively.
\end{fact}

\begin{figure}
\begin{center}
   \input{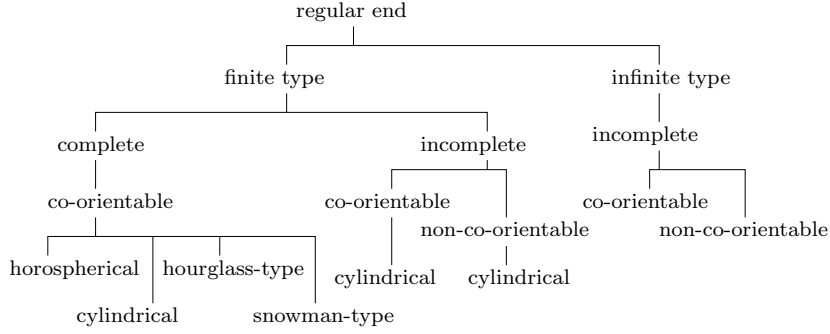}
\end{center}
\caption{Classification of weakly complete regular ends.}
\end{figure}
Let $f\colon{}M^2\to H^3$ be a weakly complete flat front with a 
regular end $p$.  We define 
\begin{equation}\label{eq:gauss-ratio}
\alpha(p) = \begin{cases}
    (dG/dG_*)(p) & 
    \text{if $\bigl | (dG/dG_*)(p) \bigr |  \le \bigl | (dG_*/dG)(p)
    \bigr |$}, \\
   (dG_*/dG)(p) &
   \text{if $\bigl | (dG/dG_*)(p) \bigr |  > \bigl | (dG_*/dG)(p)
    \bigr |$}. 
   \end{cases} 
\end{equation}
This number is called the {\it ratio of the Gauss maps\/} at the end
$p$.
The value of $\alpha$ plays an important role in criteria for the shape of 
regular ends: 
\begin{proposition}\label{fact:end-alpha}
 Let $p$ be a weakly complete regular end of a flat front $f$. 
 Then $\alpha(p)$ is contained in $[-1,1]$, 
 and is independent of rigid motions of $H^3$. 
 Moreover, $p$ is    
 \begin{enumerate}
  \item not of finite type if and only if $\alpha(p) =1$, 
  \item horospherical if and only if $\alpha(p)=0$, 
  \item of snowman-type if and only if
	$\alpha(p) >0$ and $\alpha(p)\ne 1$, 
  \item of hourglass-type if and only if
	$\alpha(p) < 0$ and $\alpha(p)\ne -1$,  
  \item cylindrical if and only if
	$\alpha(p)=-1$.  
 \end{enumerate}
 In particular, the end $p$ is complete if 
 $\alpha(p)\ne \pm 1$.
\end{proposition}
\begin{proof}
 Since the proof of Lemma~3.10 of \cite[page~165]{KUY2}
 is still valid for incomplete or non-finite-type 
 regular ends, that is, $G(p)=G_*(p)$ holds at any regular
 end $p$,
 \eqref{eq:g-moebius} implies
 $(dG/dG_*)(p)$ is invariant  under rigid motions of $H^3$.

 By definition, $\alpha$ is invariant under the operation
 in Remark~\ref{rem:dual}, so
 exchanging $G$ and $G_*$ and applying a rigid motion of $H^3$ 
 if necessary, 
 we can take a complex coordinate $z$ such that $z(p)=0$ and
 \begin{equation}\label{eq:g-gstar-alpha}
     G_* = z^m, \qquad G = a z^{m+l}\bigl(1+o(1)\bigr),
 \end{equation}
 where $m\geq 1$ and $l\geq 0$ are integers and $a\in\C\setminus\{0\}$.

 Then $\xi$ in \eqref{eq:xi} is given as
 \begin{equation}\label{eq:xi-alpha}
  \xi = c\exp\int_{z_0}^z
	      \frac{-a(m+l)z^{l-1}}{1-az^l+o(z^l)}
	      \,dz\qquad
	      (c\in \C\setminus\{0\})
 \end{equation}
 and the Hopf differential is calculated by \eqref{eq:data-g} as
 \begin{equation}\label{eq:Q-alpha}
  Q = -\frac{m(m+l)az^{l-2}\bigl(1+o(1)\bigr)}{%
     \bigl(1-az^l+o(z^l)\bigr)^2}\,dz^2.
 \end{equation}

 If $l\geq 1$,
 then $dG/dG_*=0$ at $p$, that is $\alpha(p)=0$.
 In this case, $\ord_0 Q\geq -1$ because of \eqref{eq:Q-alpha}.
 Hence $p$ is a horospherical end.
 Conversely, $\ord_0 Q\geq -1$ implies $l\geq 1$.

 Next, suppose $l=0$ and $a=1$.
 In this case, $\alpha(p)=1$ and 
 \eqref{eq:xi-alpha} implies that
 $\xi$ has an essential singularity at $z=0$.
 Then by \eqref{eq:data-g}, $\omega$ has an essential singularity
 at $0$, which implies that $p$ is not a finite type end.
 Conversely, $p$ not of finite type will similarly imply $l=0$ and
 $a=1$.

 Finally, we assume $l=0$ and $a\neq 1$.
 By the period condition in Fact~\ref{fact:g-rep} and 
 \eqref{eq:xi-alpha},
 we have $a\in \R$, 
 which implies $\alpha(p)\in\R$.
 Moreover, exchanging $G$ and $G_*$ and rechoosing $z$
 if necessary, 
 we may assume 
 \begin{equation}\label{eq:alpha-a}
   \alpha(p)=\left.\frac{dG}{dG_*}\right|_{z=p} = 
    a \in [-1,0)\cup (0,1).
 \end{equation} 
 In this case, \eqref{eq:Q-alpha} implies $\ord_0Q=-2$.
 Hence the end $0$ is cylindrical if and only if 
 $\ord_0|\omega|^2=-1$.
 Substituting \eqref{eq:xi-alpha} with $l=0$, $a=\alpha(p)$
 into the first equation of \eqref{eq:data-g}, we have
 \[ 
    \omega = z^{m\frac{1+\alpha}{1-\alpha}-1}
               \bigl(b+o(1)\bigr)\,dz
    \qquad (b\in\C\setminus\{0\}).
 \]
 Hence the end
 $0$ is cylindrical if and only if $\alpha(p)=-1$.
 Otherwise, the top term coefficient of $Q$ at $0$ 
 is obtained as $q_0=-\alpha(p) m^2 /\bigl(1-\alpha(p)\bigr)^{2}$.
 So we have the conclusion.
\end{proof}
From here on out, 
we study the ends of caustics $C_f$, which arise from 
the umbilic points and ends of $f$.  
In the former case, we call them {\em U-ends}, and in the latter case,
we call them {\em E-ends}.  We consider U-ends first:  

\begin{theorem}[Properties of U-ends]
\label{thm:u-end}
 Let $f \colon M^2 \to H^3$ be a 
 non-totally-umbilic flat front and let $p$ 
 be a point in $M^2$.  Let $C_f$ be the caustic of $f$.  
 \begin{enumerate}
  \item\label{item:u-end-1}
        If $p$ is an umbilic point of $f$, 
        then $p$ is a regular end of $C_f$ with multiplicity 
        \[
         m(C_f,p)=(\ord_p Q)/2. 
        \]
        In particular, $p$ is non-co-orientable if and only if 
        $\ord_p Q$ is odd.
        Moreover, $p$ is a cylindrical end of finite type, 
        and the singular set of $\Cf$ accumulates at $p$.  
        However, the end $p$ of $\Cf$ cannot be an end of a cylinder itself.
  \item\label{item:u-end-2} 
       Conversely, if $p$ is an end of $\Cf$, then 
        it is an umbilic point of $f$.  
 \end{enumerate}
\end{theorem}
\begin{proof}
\ref{item:u-end-1}
 Taking a rigid motion of $H^3$, if necessary, 
 we may assume $G(p)$, $G_*(p)$ are both finite.  
 Since $p \in M^2$ is not an end, $G(p)$ and $G_*(p)$ do not 
 coincide.  Thus, 
 \begin{equation}\label{g-gstar}
  \left(G - G_* \right)(p) \ne 0, \; \infty. 
 \end{equation}
 Since $p$ is an umbilic point, we have $Q(p)=0$. 
 It follows from \eqref{eq:data-g} and \eqref{g-gstar} 
 that $dG\, dG_* |_p =0$. 
 Therefore, at least one of $dG|_p$ and $dG_*|_p$ is zero. 
 We may assume that $dG|_p =0$
 (if necessary, we take the dual $\E^{\natural}$ instead of $\E$). 
 Then $dG_* |_p \ne 0$, 
 because $\G_f = (G, G_*)$ is an immersion (Theorem 2.9 of
 \cite{KUY2}).  

 Using another rigid motion of $H^3$, if necessary, we can take a local 
 coordinate $z$ centered at $p$, i.e. $z(p)=0$, such that 
 \begin{equation}\label{gstar}
  G_*(z)=z. 
 \end{equation}
 With this coordinate $z$, the hyperbolic Gauss map $G$ expands as 
 \begin{equation}\label{g}
  G(z) = a_0 + a_m z^m + a_{m+1}z^{m+1} + \cdots
     \qquad (a_0,a_m\neq 0),
 \end{equation} 
 where $m \ge 2$. 
 We remark that 
 $ \ord_p Q = m-1$, 
 and $dG / dG_*$ is  computed as 
 \begin{equation}\label{alpha}
  dG / dG_* = m a_m z^{m-1} h(z)\qquad (h(0)=1), 
 \end{equation}
 where $h(z)$ is holomorphic in $z$.  
 In particular, $(dG / dG_*)(0) =0$. 
 It follows from \eqref{hGm_c}, \eqref{gstar}, \eqref{g}, 
 \eqref{alpha} that 
 $\Gc(p) = \Gc_*(p) \left( = G(p) \right)$, so $p$ is an end of $\Cf$. 
 Then \eqref{hGm_c}, \eqref{gstar}, \eqref{g} imply 
 $\Gc$ and $\Gc_*$ are meromorphic at $p$ (on the double cover
 of a neighborhood of $p$), so $p$ is a regular 
 end of $\Cf$.  In particular, 
 substituting \eqref{gstar}, \eqref{g}, \eqref{alpha} into 
 \eqref{hGm_c}, we have 
 \begin{align*}
  \Gc & = a_0-a_0 \sqrt{m a_m}z^{(m-1)/2}h_1(z)+
           o\left(z^{(m-1)/2}\right),\\
  \Gc_* &=  a_0+a_0 \sqrt{m a_m}z^{(m-1)/2}h_1(z)+
           o\left(z^{(m-1)/2}\right),
 \end{align*}
 where $h_1(z)$ is a holomorphic function in $z$ such that $h_1{}^2=h$,
 with $o(z^{(m-1)/2})$ denoting higher order terms.  
 The multiplicity of the end $p$ is 
 \begin{equation*}
  m(\Cf,p) = \tfrac{1}{2}(m-1) = \tfrac{1}{2}\ord_p Q. 
 \end{equation*}

 It follows from \eqref{omega_c}, \eqref{alpha} 
 that 
 \begin{equation}\label{omegac}
  4 \omegac = \frac{2 ( \beta+1 )^2}{G-G_*}dz
  -(m-1)\frac{dz}{z}-\frac{h'}{h}dz 
       = 
         \frac{1}{z}
          \bigl( 1-m + o(1)\bigr) dz.
 \end{equation}
 Similarly, by \eqref{theta_c}, \eqref{alpha}, 
 \begin{equation}\label{thetac}
  4 \thetac =
        \frac{1}{z}
	\bigl(1-m + o(1) \bigr) dz . 
 \end{equation}
 Hence, $\ord_p |\omegac|^2 = \ord_p |\thetac|^2 = -1$. 
 This implies that $p$ is a weakly complete cylindrical end of finite type. 
 Moreover, 
 \[
  \rhoc = \frac{\thetac}{\omegac} 
   =\frac{m-1+o(1)}{m-1+o(1)}, 
 \]
 in particular,  $\rhoc(p)=1$.
 So the singular set $\{ |\rhoc| = 1\}$ accumulates at $p$. 
 
 Finally, we show that $p$ is not an end of 
 a hyperbolic cylinder. 
 Suppose, by way of contradiction, that $p$ is an end of revolution, so 
 $\thetac = k \omegac$ for some $k \in \C$, 
   in a neighborhood $U_p$ of $p$.  
 Comparing the coefficients of $1/z$ in \eqref{omegac} and
 \eqref{thetac}, 
 $k$ is necessarily $+1$ and 
 $\thetac = \omegac$. 
 Then \eqref{omega_c}, \eqref{theta_c} give 
 $dG/dG_* = 0$ on $U_p$, contradicting \eqref{alpha}. 

\ref{item:u-end-2}
 Suppose now $p \in M^2$ is an end of $\Cf$, that is, 
 $\Gc(p) = \Gc_*(p)$.   
 Without loss of generality, we may assume $G(p) \ne \infty$, 
 $G_*(p) \ne \infty$ and $G(p) \ne G_*(p)$. 
 Then \eqref{hGm_c} gives $(dG/dG_*)(p) = 0$ or $\infty$, so 
 one of $dG$, $dG_*$ is zero at $p$, and so $p$ is umbilic. 
\end{proof}

\begin{remark}\label{remarkaddedbywayne}
The claim in Theorem \ref{thm:u-end} that $p$ is a non-co-orientable end of 
$\Cf$ if and only if $\ord_p Q$ is odd follows in this way: 
The end $p$ is non-co-orientable if and only if the deck 
transformation associated to a once-wrapped loop about $p$ switches $\omegac$ and 
$\thetac$ with respect to a local adjusted lift, by Corollary \ref{cor:double-lift}.  
Then non-co-orientability is equivalent to $\ord_p Q$ being odd, by 
\eqref{omega_c-theta_c}.  
\end{remark}

We shall next consider E-ends.  
To avoid confusion, we denote by $(f;p)$ the end $p$ of $f$,  
and by $(\Cf;p)$ the end $p$ of $\Cf$.  

Suppose that the multiplicity of the end $(f;p)$ is $m$, i.e.,  
$m(f,p)=m \in \Z_{+}$. 
Without loss of generality, 
we may assume that $G(p)=G_*(p)=0$, and that 
$r_p(G_*)=m$, $ r_p(G)=m+k$, where $k$ is a non-negative integer.  
There exists a local coordinate $z$ centered at $p$ 
such that 
\begin{equation}\label{eq:g-star-end}
 G_*(z) = z^m .
\end{equation}
With this coordinate $z$, $G$ expands, with $h(0) = a_{m+k} \neq 0$, as 
\begin{equation}\label{eq:g-end} 
 G(z) = a_{m+k} z^{m+k} + a_{m+k+1}z^{m+k+1} + \cdots = z^{m+k} h(z) \; . 
\end{equation}

\subsection*{The case of ``$k > 0$'' or ``$k=0$ with $a_m \ne 1$'' in 
\eqref{eq:g-end}}  
By \eqref{eq:data-g}, we have 
\[
 Q = -\frac{m z^{k-2}h_1(z)}{(z^{k}h(z)-1)^2} \, dz^2 ,
\]
where $h_1(z) = (m+k) h(z)+zh'(z)$. 
Since
$h_1(0) = (m+k) h(0) = (m+k)a_{m+k} \neq 0$, 
\begin{equation}\label{ord_pQ}
 \ord_p Q = k-2 \ge -2.
\end{equation}

We can also compute that
\begin{equation}\label{alpha-2}
 \frac{dG}{dG_*} = \frac{1}{m}z^{k}h_1(z).   
\end{equation}
The equations \eqref{hGm_c}, \eqref{alpha-2} yield  
\[
 \Gc(z) = \frac{%
       z^{(2m+k)/2}%
       \left\{ \sqrt{m} \, z^{k/2} h + \sqrt{h_1} \right\} }
       {\sqrt{m} + z^{{k/2}}\sqrt{h_1}}, 
    \quad
 \Gc_*(z) = \frac{%
       z^{{(2m+k)/2}}%
      \left\{ \sqrt{m} \, z^{{k/2}} h - \sqrt{h_1} \right\} }
       {\sqrt{m} - z^{{k/2}}\sqrt{h_1}}. 
\]
Therefore, 
$r_p(\Gc)=r_p(\Gc_*)= m+(k/2)$.
It follows from \eqref{ord_pQ} that
\[
 m(\Cf,p)= m+(k/2) = (1/2) \ord_p Q +m+1. 
\]
The equations 
\eqref{omega_c}, \eqref{theta_c} and \eqref{alpha-2} yield
\begin{align}
 4 \omegac &= \left( \frac{1}{z} \left\{ 
   \frac{2(z^{k/2}\sqrt{h_1}+\sqrt{m})^2}{z^{k}h-1}-k
   \right\}
 -\frac{h_1'}{h_1} \right)dz, \label{omega_c_2} \\ 
 4 \thetac &= \left( \frac{1}{z} \left\{ 
 \frac{2(z^{k/2}\sqrt{h_1}-\sqrt{m})^2}{z^{k}h-1}-k
 \right\}
 -\frac{h_1'}{h_1} \right) dz.  \label{theta_c_2}
\end{align}
These imply that
$ \ord_p |\omegac|^2 = \ord_p |\thetac|^2 = -1$, 
hence $(\Cf;p)$ is a cylindrical end of finite type. 
Moreover, 
by \eqref{omega_c_2}, \eqref{theta_c_2}, we can prove 
\begin{equation*}
 \lim_{z \to 0} \frac{\thetac}{\omegac} = 
 \begin{cases}
  1 & \text{if } k > 0,  \\
  \left(\dfrac{\sqrt{a_m}-1}{\sqrt{a_m}+1}\right)^2 & 
  \text{if } k = 0, \ a_m \ne 1. 
 \end{cases}
\end{equation*} 
Thus, the singularities of $\Cf$ accumulate at $p$ (i.e., 
$\lim_{z \to 0} \left| \thetac / \omegac  \right| =1$)  
if and only if $k > 0$ or $k=0$ with negative real number $a_m$. 

\subsection*{The case of ``$k=0$ with $a_m = 1$'' in \eqref{eq:g-end}}
Now $G$ is written as
\begin{equation*}
 G(z) = z^m + a_{l}z^{l} + a_{l+1}z^{l+1} + \cdots 
      = z^m + a_l z^l +o(z^{l}) \; , 
\end{equation*}
where $l > m$, $a_{l} \ne 0$, and 
$o(\cdot)$ denotes higher order terms in $z$.  
One easily gets 
\[
 \ord_p Q = 2(m-1)-2 l = 2(m-l-1) \le -4. 
\]
In particular, $\ord_p Q$ is an even number.  
On the other hand,
\begin{equation}\label{eq:end-g-ratio}
 \frac{dG}{dG_*} = \frac{m z^{m-1}+ la_lz^{l-1}+\cdots}{mz^{m-1}}
                 = 1+\frac{l}{m}a_lz^{l-m}+o(z^{l-m}),
\end{equation}
\begin{equation}\label{eq:end-beta}
  \beta = \sqrt{\frac{dG}{dG_*}} =
          1+\frac{l}{2m}a_l z^{l-m} + o(z^{l-m}).
\end{equation}
Then \eqref{hGm_c}, \eqref{eq:end-beta} yield 
\[ \Gc(z) = z^m \left(
                1+o(1)
               \right),\qquad
 \Gc_*(z) = z^m \left(
                \frac{l-2m}{l} + o(1) 
               \right) \; . \]
Hence, 
regardless of whether $l = 2m$ or $l \ne 2m$, it follows that
$\min \{ r_p(\Gc), \,  r_p(\Gc_*)\} =m$.  
Therefore we have, for any $l$, that 
\[
 m(\Cf,p)=m \; . 
\]

Next, we investigate $\omegac$, $\thetac$ around $p$.  
It follows from \eqref{eq:end-g-ratio} that 
\[
 d \log \left(\frac{dG}{dG_*}\right) = 
    \left(\frac{l(l-m)}{m}a_l z^{l-m-1}+ o(z^{l-m-1})\right)\,dz.
\]
Then by \eqref{omega_c}, \eqref{theta_c} and \eqref{eq:end-beta},
we have
\[ 
  4 \thetac=z^{l-m-1}\left(\frac{l(2m-l)}{2 m}a_l+o(1)\right) dz , \qquad
  4 \omegac = z^{m-l-1}\left(\frac{8m}{a_l}+o(1)\right) dz.
\]
Hence
\[
   \ord_p |\omegac|^2=m-l-1\quad\text{and}\quad
   \ord_p |\thetac|^2
     \begin{cases} 
      =l-m-1 \quad & \text{if $l\neq 2m$}\\
      >l-m-1=m-1   \quad & \text{if $l=2m$}
     \end{cases},
\]
so $\ord_p |\thetac|^2\geq 0$, $\ord_p|\omegac|^2\leq -2$.  
Hence $(C_f;p)$ is a non-cylindrical end of finite type.  

\subsection*{Summary of the above argument} Under the situation above, 
 \begin{itemize}
  \item Assume that $k>0$, or $k=0$ with $a_m \ne 1$. Then 
        \begin{enumerate}[\rm (i)]
         \item $\ord_p Q = k -2 \ge -2$, 
         \item $m(\Cf,p) = \frac{1}{2}\ord_p Q + m+1$, 
         \item $(\Cf;p)$ is a cylindrical end of finite type, 
         \item the singularities of $\Cf$ accumulate at $p$ if and only if 
               $k>0$, or $k=0$ with $a_m < 0$.  
        \end{enumerate}
  \item Assume that $k=0$ with $a_m = 1$. Then 
  \begin{enumerate}[\rm (i)]
   \item $\ord_p Q = 2(m-l-1) \le -4$, 
   \item $m(\Cf,p) = m$,
   \item $(\Cf;p)$ is a non-cylindrical end of finite type, 
   \item the singularities of $\Cf$ do not accumulate at $p$.  
  \end{enumerate}
 \end{itemize}

Using Fact~\ref{fact:end-alpha} and Remark 
\ref{remarkaddedbywayne}, we can restate these conclusions as: 

\begin{theorem}[Properties of E-ends]\label{thm:e-end}
 Let $f\colon M^2 \to H^3$ be a non-totally-umbilic
 weakly complete flat front, with a regular end 
 $p$ {\rm (}see Definition \ref{def:puncture}{\rm )}. 
 Then $p$ is also a weakly complete regular end of $C_f$.  Moreover: 
 \begin{itemize}
  \item If $\ord_p Q \ge -2$, then the end $(\Cf;p)$ has multiplicity 
        $m(\Cf,p) = \frac{1}{2}\ord_p Q + m(f,p) + 1$, and 
        $(\Cf;p)$ is non-co-orientable if and only if $\ord_p Q$ is odd. 
        Moreover, $(\Cf;p)$ is a cylindrical end of finite type. 
        The singularities of $\Cf$ do not accumulate at $p$ if and only if 
        $(f;p)$ is of snowman-type.  
  \item If $\ord_p Q < -2$, then 
        $\ord_p Q \leq -4$ and is necessarily an even integer, and 
        the end $(\Cf;p)$ has multiplicity 
        $m(\Cf,p) = m(f,p)$. In particular, $(\Cf;p)$ is co-orientable. 
        Moreover, $(\Cf;p)$ is a non-cylindrical end of finite type. 
 \end{itemize}
\end{theorem}
Finally, we prove Theorem~\ref{thm:C} in the introduction.
\begin{proof}[Proof of Theorem~\ref{thm:C}]
 (2) follows from (1) by 
 Theorems~\ref{thm:u-end}, \ref{thm:e-end}. 
 We now assume (2) and prove (1).  By Proposition 
 \ref{inserted-at-last-minute-proposition}, 
 the domain of $\Cf$ is biholomorphic 
 to a compact Riemann surface minus finitely many points, so the same
 holds for  $f$ as well.  

 Let $p$ be an arbitrary end of $\Cf$. 
 By Corollary \ref{cor:double-lift} and Proposition \ref{prop:finite->} and 
 assumption (2), $\Cf$ lifts to a double cover on which its canonical
 $1$-forms $\omegac$ and $\thetac$ have finite order at $p$.  
 By \eqref{omega_c-theta_c}, 
 $Q = -(\omegac - \thetac)^2/4$, and $\ord_p Q$ is finite.
 Since $\omegac$ has finite order and $\Gc$ is meromorphic, each
 component of $\Ec$ 
 has finite order at $p$, by \eqref{eq:sm-analogue2} for $\Cf$.  
 By \eqref{eq:Ec}, the components $(\Ec)_{ij}$ of $\Ec$ satisfy 
 \begin{align*}
  \bigl( (\Ec)_{21} + (\Ec)_{22} \bigr)^4 &= 
   \left( \frac{\sqrt{2} i}
    {\sqrt{(G - G_*) \beta}} \right)^4 = -
  \frac{4Q}{dG^2}\qquad\text{and}\\
  \bigl( (\Ec)_{21} - (\Ec)_{22} \bigr)^4
   &= 
   \left( \frac{\sqrt{2} i \sqrt{\beta}}
    {\sqrt{G - G_*}} \right)^4 = -
   \frac{4Q}{dG_*^2}.
 \end{align*}
 Therefore 
 $dG$ and $dG_*$ have finite orders at $p$, so 
 $G$ and $G_*$ do as well.  
 Hence $f$ has a regular end at $p$.  
 
 Next we shall show that $f$ is weakly complete (at any E-end $p$). 
 Without loss of generality, we may assume $G(p) = G_*(p)=0$. 
 If $\ord_p Q\le -2$, then 
 \[
    ds^2_{1,1}=|\omega|^2+|\theta|^2
              =(|\hat\omega|^2+|\hat\theta|^2)|dz|^2
               \ge 2|\hat\omega\hat\theta||dz|^2=2|\hat Q||dz|^2, 
 \]
 where $\omega=\hat\omega\,dz$, $\theta=\hat\theta\,dz$ and
 $Q=\hat Q \,dz^2$.
 So obviously $f$ is weakly complete at $p$. 
 Let us consider the case $\ord_p Q\ge -1$. 
 Then there is a local coordinate $z$ giving 
 \begin{equation}\label{eq:GGstar_final}
  G_*=z^m\quad\text{and}\quad
  G=z^{m+k}\bigl(a+o(1)\bigr), \qquad 
   (a\neq 0,k\ge 1).
 \end{equation}
 Therefore, 
 ${dG}/{(G-G_*)}$
 has order $k-1$ at $z=0$, and so is holomorphic at $z=0$.  Then 
  $\xi$ in \eqref{eq:xi} 
 is a holomorphic function which does not vanish at $z=0$.
 This implies that $f$ is weakly complete at $z=0$, 
 because of \eqref{eq:data-g} and \eqref{eq:GGstar_final}. 
\end{proof}

\end{document}